\documentclass[12pt]{article}
\usepackage[myheadings]{fullpage}
\usepackage{fancyhdr}
\usepackage{amsmath}
\usepackage{amsthm}
\usepackage{amssymb}
\usepackage{lastpage}
\usepackage{hyperref}
\usepackage{graphicx}
\usepackage{datetime}
\usepackage{pdfpages}
\usepackage{caption}
\usepackage{booktabs}
\usepackage{float}
\usepackage{graphicx}
\usepackage{subcaption}
\usepackage{tikz}
\usepackage{natbib}
\usepackage{url}
\usepackage{authblk}
\usetikzlibrary{calc, arrows.meta}

\usepackage{xcolor}
\hypersetup{colorlinks, citecolor=blue}

\newcommand{\abs}[1]{\left\vert#1\right\vert} %$\abs x$
\newcommand{\set}[2]{\left \{ \left. #1 \,\right| #2 \right\}  } %$\set{x}{x\in\QQ}$
\newcommand{\injr}[1]{\text{inj}#1}
\newcommand{\norm}[1]{\left\lVert#1\right\rVert} % norm/length notation
\newcommand{\inpr}[2]{\left\langle #1, #2 \right\rangle} % inner product notation

\newcommand{\var}{\operatorname{Var}}

\newcommand{\pd}[2]{\frac{\partial #1}{\partial #2}}
\providecommand{\keywords}[1]
{
  \small	
  \textbf{Keywords:} #1
}

\newtheorem{theorem}{\underline{Theorem}}[section]
\newtheorem{proposition}[theorem]{\underline{Proposition}}
\newtheorem{lemma}[theorem]{\underline{Lemma}}

\theoremstyle{definition}
\newtheorem{definition}[theorem]{\underline{Definition}}

\newtheorem{example}[theorem]{\underline{Example}}

\title{Autoregressive Processes on Riemannian Manifolds}

\author[1,2]{Meshal Abuqrais\thanks{Emails: \href{mailto:meshal.abuqrais@kcl.ac.uk}{\texttt{meshal.abuqrais@kcl.ac.uk}}, \href{mailto:M.abuqrais@ku.edu.kw}{\texttt{m.abuqrais@ku.edu.kw}}}}
\author[1]{Davide Pigoli\thanks{Email: \href{mailto:davide.pigoli@kcl.ac.uk}{\texttt{davide.pigoli@kcl.ac.uk}}}}

\affil[1]{King's College London}
\affil[2]{Kuwait University}

\begin{document}

\maketitle
%\tableofcontents

\begin{abstract}
This paper introduces a Riemannian autoregressive (R-AR) model of order one for manifold-valued time series. 
The model is specified through an autoregressive process in the tangent space at a reference point, which is mapped to the manifold via the exponential map. 
It is characterised by two parameters: a reference point $\mu$, representing the intrinsic central tendency through the Fréchet mean, and an autoregressive parameter $\phi$, governing the dependence structure and stationarity properties.
When these parameters are unknown, their estimation introduces geometric and probabilistic challenges arising from the intrinsic estimation of $\mu$ and from estimating $\phi$ as the base point varies over the manifold and the associated tangent spaces change.
The estimation theory is developed using a strong law of large numbers for the sample Fréchet mean set of ergodic Markov chains, together with ergodic arguments for the strong consistency of the autoregressive parameter estimator.
The framework is validated through numerical simulations in the hyperbolic plane and an application to aerosol size distributions on the Fisher-Rao manifold, demonstrating how the proposed model can characterise mean-reverting dynamics in nonlinear geometries.
\end{abstract}


\section{Introduction}
Classical time series models, such as autoregressive stochastic processes, have been widely applied across fields, from economics \cite{cipra2020time} to medicine \cite{jelinek2017ecg}, and remain of significant theoretical interest. 
Typically, such models are formulated in linear spaces, such as Euclidean or Hilbert spaces, which provide a linear structure that facilitates their analysis. 
In recent years, however, there has been growing interest in analysing data that take values in non-Euclidean spaces, where the underlying geometry fundamentally influences statistical methodology. 
For example, the work \cite{ding2025manifold} studies time series of symmetric positive-definite matrices for EEG data, \cite{vemulapalli2016rolling} models human skeletal motion as points on a manifold of rotations, and \cite{nava2025ridge} analyses spherical data for hurricane tracking. 
In such settings, traditional Euclidean-based models become inadequate, necessitating statistical frameworks that respect the intrinsic geometry of the data.

One of the primary obstacles to analysing data in more general spaces, particularly manifold-valued data, is the lack of an intrinsic algebraic or vector space structure, which prevents the direct application of classical linear methods.
Nevertheless, variational characterisations of basic statistical quantities, such as the mean and variance, remain applicable because of their compatibility with manifold or metric-space structures. 
The most prominent example is the Fréchet mean \citep{Frechet1948elements}, a generalisation of the classical mean defined through a variational principle that is well suited to nonlinear spaces.

Stochastic processes on manifolds have been extensively studied, particularly in the form of Brownian motion and diffusion processes; see, for example, \cite{hsu2002stochastic} and \cite{ikeda2014stochastic}. 
These continuous-time models are well suited for capturing local stochastic behaviour, but they do not provide a natural framework for modelling discrete-time dependence, long-term memory, or mean-reverting dynamics, all of which are central to many applications in time series analysis.
Autoregressive models address these needs directly but, in the non-Euclidean setting, require a new formulation to define central tendencies and assess stability. Motivated by this gap, we develop a framework for manifold-valued autoregressive processes as a generalisation of their Euclidean discrete-time counterparts.

This paper investigates a class of autoregressive stochastic processes defined on Riemannian manifolds that is analogous to first-order autoregressive models.
While extensions of classical autoregressive models to non-Euclidean settings have previously been proposed, much of the existing work focuses on specific types of manifolds, such as Lie groups, spheres, or matrix manifolds, where additional structure, such as group operations or matrix algebra, can be exploited. 
For instance, \cite{fiori2014auto} develops an ARMA-type model on matrix manifolds by exploiting their algebraic properties. 
In \cite{zhu2024spherical}, a model for autoregressive processes of sphere-valued random variables is constructed by exploiting spherical geometry.

Our work instead builds on the recent approach of \cite{bulte2024autoregressive}, which studies autoregressive processes in Hadamard geodesic spaces from the perspective of random dynamical systems.
In contrast, we develop this perspective on Riemannian manifolds, thereby accommodating richer geometrical structures.
Our approach to defining the model relies on the tangent space at a reference point, which plays the role of the Frechet mean, and on the global injectivity radius.
For a fixed reference point $\mu$, the process corresponds to an autoregressive process in $T_\mu M$ mapped to the manifold through the exponential map.
This allows us to define the model locally even when the space has positive curvature.
Furthermore, we view the model as a Markov chain on a Riemannian manifold induced by a random dynamical system and analyse it using Markov chain theory.
The geometric aspects become more pronounced in the estimation problem, where the reference point is unknown and the estimation of $\phi$ must account for the variation of the associated tangent spaces as the base point changes.
By integrating tools from stochastic processes and Riemannian geometry, we define and analyse manifold-valued autoregressive processes. 
The proposed model is characterised by two parameters, which can be interpreted as a Fréchet mean term representing central tendency and a contraction parameter governing stability and mean reversion.

The main contributions of this article are therefore twofold.
First, we introduce a first-order Riemannian autoregressive model as a generalisation of classical autoregressive processes, formulate it as a random dynamical system on a manifold, and propose suitable estimators for the model parameters.
Second, we use a strong law of large numbers for the sample Fr\'echet mean of an ergodic Markov chain as a basis for establishing the strong consistency of the proposed parameter estimators.

The paper is organised as follows. 
Section 2 provides the necessary background on differential geometry, probability theory, and the Fréchet function. 
Section 3 formally defines the autoregressive process on Riemannian manifolds and establishes its basic properties. 
Section 4 establishes the asymptotic properties of the proposed estimators of the model parameters.
Section 5 presents numerical simulations and an application to a real dataset in hyperbolic spaces.
Section 6 concludes the paper and discusses directions for future work.

\section{Geometric and Probabilistic Background}

\subsection{Riemannian manifolds}
A $d$-dimensional smooth manifold, referred to simply as a manifold, is a space that locally can be modelled as a $d$-dimensional Euclidean space, allowing the tools of calculus to be extended beyond the classical Euclidean setting.
At each point, the manifold can be approximated by a vector space called the tangent space at that point. 
A Riemannian metric assigns an inner product to each tangent space, thereby making it possible to measure lengths, angles, and distances.
A manifold equipped with such a metric is called a Riemannian manifold.
In this subsection, we introduce the necessary background on manifolds that will be used throughout the article. For a comprehensive treatment of Riemannian geometry, we refer the reader to \cite{Rie-Lee,jost2008riemannian}.
We always assume that the manifolds are finite-dimensional.

Formally, a $d$-dimensional Riemannian manifold is a $d$-dimensional smooth manifold $M$ equipped with a smoothly varying inner product $\inpr{\cdot}{\cdot}_p$ on each tangent space $T_pM$, for $p \in M$. 
Such a family of inner products is called a Riemannian metric on $M$.
We denote the induced norm on the tangent space $T_pM$ by $\norm{\cdot}_p$ or $\norm{\cdot}$ if the tangent space is clear from the context.
The tangent bundle of $M$ is denoted by $TM$ and defined as the disjoint union of all tangent spaces of $M$.
The tangent bundle has a natural smooth structure that is constructed from $M$ such that the projection $\pi:TM\to M, (v;p)\mapsto \pi(v;p)=p$ is smooth for all $p$.

A smooth vector field, which we simply call a vector field, is a smooth map $V:M\to TM$ such that $\pi\circ V=\mathrm{id}_M$.
The set of all smooth vector fields on $M$ is denoted by $\mathfrak{X}(M)$.
If $I$ is an interval in $\mathbb{R}$, a smooth curve, or simply a curve, $\gamma$ is a smooth map from $I$ to $M$.
A connection on a smooth manifold $M$ is a map 
\begin{equation*}
\nabla : \mathfrak{X}(M) \times \mathfrak{X}(M) \to \mathfrak{X}(M),
\end{equation*}
that assigns to two vector fields $X$ and $Y$ another vector field $\nabla_X Y$, called the covariant derivative of $Y$ in the direction of $X$. 
This operator is linear in $X$ and satisfies the product rule
$\nabla_X(fY) = (Xf)Y + f \nabla_X Y,$
for any smooth function $f$ on $M$.
The connection is called symmetric if the Lie bracket of any two vector fields $X$ and $Y$ on $M$ is given by $[X, Y] = \nabla_X Y - \nabla_Y X$.
In principle, connections are independent of the Riemannian metric.
A cornerstone theorem in Riemannian geometry states that for a Riemannian metric $\inpr{\cdot}{\cdot}$ on $M$, there exists a unique symmetric connection that is compatible with the metric, meaning
\begin{equation*}
Z \inpr{X}{Y} =  \inpr{\nabla_Z X}{Y} +\inpr{X}{\nabla_Z Y},
\end{equation*}
for all vector fields $X, Y, Z$ on $M$.
This unique connection is called the Levi-Civita connection of $(M,\inpr{\cdot}{\cdot})$, and we assume that our Riemannian manifolds are equipped with their Levi-Civita connections.
A smooth curve $\gamma$ in $M$ is said to be a geodesic if
\begin{equation*}
\nabla_{\gamma'}\gamma'=0.
\end{equation*}
The Riemannian metric induces a natural distance function on $M$ through the lengths of paths.
Assume $M$ is a connected manifold.
The length of a piecewise smooth curve on $[0,1]$ is given by 
\begin{equation*}
L(\gamma)=\int_{0}^1\norm{\gamma'(t)}_{\gamma(t)}dt.
\end{equation*}
The distance between $p$ and $q$ is $d(p,q)=\inf_{\gamma}{L(\gamma)}$, where the infimum is taken over all piecewise smooth paths from $p$ to $q$.
This turns $M$ into a metric space, and whenever we speak of a Riemannian manifold as a metric space, we are referring to this distance function.

We now turn to the construction of the exponential map on $M$.
The exponential map provides natural local coordinates for $M$.
Fix a point $p$ in $M$ and let $v\in T_pM$ be any tangent vector to $M$ at $p$.
Then there exists a unique maximal geodesic $\gamma_v(t;p)$ defined on an interval in $\mathbb{R}$ containing zero such that $\gamma_v(0;p)=p$ and $\gamma'_v(0;p)=v$.
Furthermore, one can rescale the interval of the maximal geodesic to include $1$.
We can then define the exponential map as
\begin{equation*}
\exp_{p}(v)=\gamma_v(1;p),
\end{equation*}
provided that $\gamma_v(1;p)$ is defined.
Note that in general, the exponential is not defined on the entire tangent bundle $TM$, but rather on a subset of it.
Specifically, for every $p$, the exponential map is defined on a neighbourhood of $0\in T_pM$.
The exponential map restricts to a local diffeomorphism near $0\in T_pM$.
The maximal radius of an open ball centred at $0$ in $T_pM$ on which the exponential is a diffeomorphism is called the injectivity radius of $M$ at $p$, written as $\injr{p}$.
The global version is defined by $\injr{M}=\inf_{p\in M}{\injr{p}}$.
Within the injectivity radius, the exponential map is invertible, and its inverse is called the logarithm map, written as $\log_p(\cdot)$.
Within the injectivity radius $r$, the exponential map is a radial isometry.
Specifically, for all $q\in B_{p}(r)$ where $r<\injr{p}$
\begin{equation*}
d(p,q)=\norm{\log_{p}q}.
\end{equation*}
A fundamental theorem in Riemannian geometry, the Hopf–Rinow theorem, establishes the equivalence between metric completeness of $M$ and the exponential map being defined on the entire tangent bundle.
In fact, these properties are also equivalent to $M$ possessing the Heine-Borel property: the compact subsets of $M$ are exactly those that are closed and bounded, \cite{jost2008riemannian}.

Finally, we turn to the curvature of the Riemannian manifold. 
For any $X, Y \in \mathfrak{X}(M)$, we define the curvature operator $R(X,Y): \mathfrak{X}(M) \to \mathfrak{X}(M)$ by
\begin{equation*}
R(X,Y)Z = \nabla_X \nabla_Y Z - \nabla_Y \nabla_X Z - \nabla_{[X,Y]} Z.
\end{equation*}
At any point $p \in M$, the sectional curvature of $M$ at $p$ for the plane spanned by $X, Y \in T_pM$ is defined as
\begin{equation*}
K_p(X,Y) = \frac{\inpr{R(X,Y)Y}{X}}{\norm{X}^2 \norm{Y}^2 - \inpr{X}{Y}^2}.
\end{equation*}
In general, $K$ depends on the point $p$ and the plane determined by $X,Y$. 
We say $M$ has curvature bounded above by $K_U$ if for all $p\in M$ and $X,Y\in T_{p}M$ $
K_{p}(X,Y)\leq K_U$.
A similar notion for bounded from below.
Also, $M$ has constant curvature $K$ if $K$ is the same for all $p \in M$ and all $X,Y \in T_{p}M$.
Typical examples of such manifolds are the space forms: Euclidean space $\mathbb{R}^d$ ($K = 0$), the unit sphere $\mathbb{S}^d$ with the round metric ($K = 1$), and hyperbolic space $\mathbb{H}^d$ with the canonical hyperbolic metric ($K = -1$).
We will discuss the hyperbolic space in an example later in the article.

\subsection{Fréchet function and Markov chains in metric spaces}
In this subsection, we introduce the Fréchet function and Markov chains in metric spaces.
For a comprehensive treatment of the Fréchet function of random objects and Markov chains, we refer the reader to \cite{patrangenaru2016nonparametric} and \cite{benaim2022markov}, respectively.
Throughout this section, we assume that $(M,d)$ is a Polish metric space, i.e., complete and separable, equipped with its Borel $\sigma$-algebra $\mathcal{B}(M)$.
Also, we fix a filtered probability space $(\Omega,\mathcal{A},(\mathcal{A}_n)_{n\in\mathbb{N}_0},P)$.

In the absence of linear structure, classical notions such as the mean must be defined differently.
To study distributional properties of random objects in general metric spaces, we introduce the Fréchet function.

We say that $X$ is a random object, or an $M$-valued random variable, if it is a measurable map from $(\Omega,\mathcal{A})$ to $(M,\mathcal{B}(M))$.
The Fréchet $r$-th moment function of $X$ is defined by
\begin{equation*}
\mathcal{F}_{r,X}(p) = \int_{M} d^r(x, p) P_X(dx) = E[d^r(X, p)],
\end{equation*}
where $r > 0$ and $P_X$ is the law of $X$.
Because we are primarily interested in the case $r = 2$, we will write $\mathcal{F}_X(p)$ or simply $\mathcal{F}(p)$ when $X$ is clear from context.

Due to the lack of algebraic structure in general metric spaces, the classical definition of the mean does not extend directly.
However, the variational characterisation of the mean does extend.
We define the Fréchet mean set of an $M$-valued random object $X$ as
\begin{equation*}
\mu_X = \operatorname{arg\,min}_{p \in M} \mathcal{F}(p).
\end{equation*}

Unlike in linear spaces, the Fréchet mean may not be unique; thus we speak of the Fréchet mean set.
Nevertheless, there are cases in which the Fréchet mean is unique.
For example, in complete Riemannian manifolds, under certain regularity conditions on the curvature and the support of the random object, it has been shown in \cite{afsari2011riemannian} that the Fréchet mean is unique.
In \cite{sturm2003probability}, global uniqueness is established for metric spaces with nonpositive curvature in the sense of Alexandrov.
This implies global uniqueness of the Fréchet mean in Hadamard manifolds, i.e., complete, simply connected Riemannian manifolds with nonpositive sectional curvature.
In the literature, particularly for manifold-valued random variables, local minimisers of the Fréchet function are referred to as Karcher means; see, for example, \cite{pennec2006intrinsic}.

A natural question concerns the estimation of the Fréchet mean set.
Given $X_1,\ldots,X_n$, the natural estimator is
\begin{equation*}
\mu_n = \operatorname*{arg\,min}_{p \in M} \frac{1}{n} \sum_{k=1}^n d^2(X_k,p),
\end{equation*}
known as the sample Fréchet mean set, or the empirical Fréchet mean set.
It was shown in \cite{bhattacharya2003large} that, under the assumption that $X_1,\ldots,X_n$ are independent and identically distributed, together with suitable assumptions on the space, the sample Fréchet mean set is strongly consistent under an appropriate notion of set convergence.
More recently, \cite{lee2025kolmogorov} established a Kolmogorov-Feller-type weak law of large numbers for sample Fréchet means on non-compact symmetric spaces under independence assumptions, including the independent non-identically distributed setting.
More general ergodic results for dependent sequences have recently been obtained in \cite{jaffe2024fr}.
In the setting considered here, we follow the approach in \cite{bhattacharya2003large} and adapt its consistency argument from independent observations to ergodic Markov chains.
This formulation will also serve as the basis for the subsequent strong consistency results of the proposed parameter estimators.
A Markov kernel $\mathcal{K}$ on $(M,\mathcal{B}(M))$ is a map
\begin{equation*}
\mathcal{K}\colon M\times \mathcal{B}(M)\to [0,1],
\end{equation*}
such that
\begin{enumerate}
\item For all $B\in\mathcal{B}(M)$, the map $x \mapsto \mathcal{K}(x,B)$ is measurable.
\item For all $x\in M$, the map $B \mapsto \mathcal{K}(x,B)$ is a probability measure on $\mathcal{B}(M)$.
\end{enumerate}

Let $\mathcal{K}$ be a Markov kernel on $M$.
An adapted stochastic process $(X_n)_{n\geq 0}$ of $M$-valued random variables on $(\Omega,\mathcal{A},(\mathcal{A}_n)_{n\geq 0},P)$ is called a Markov chain on $M$ with kernel $\mathcal{K}$ if for all $n\in\mathbb{N}_0$ and $B\in\mathcal{B}(M)$
\begin{equation*}
P(X_{n+1}\in B \mid \mathcal{A}_n) = \mathcal{K}(X_n,B).
\end{equation*}

A Markov kernel naturally induces a map on probability measures on $M$, $\mathsf{Prob}(M)$, by
\begin{equation*}
\mu \mapsto \mu\mathcal{K}(\cdot)=\int_{M}\mathcal{K}(x,\cdot)\,d\mu(x),
\end{equation*}
which defines a probability measure on $M$.
If there exists a probability measure $\lambda\in\mathsf{Prob}(M)$ such that $\lambda\mathcal{K}=\lambda$, then $\lambda$ is called an invariant probability measure for $\mathcal{K}$.

Any Markov chain admits a canonical realisation on the path space $M^{\mathbb{N}_0}$.
Denote by $\mathcal{B}(M)^{\otimes \mathbb{N}_0}$ the product $\sigma$-algebra.
Suppose that $\lambda_0$ is the probability law of $X_0$.
There exists a unique probability measure $\mathbb{P}_{\lambda_0}$ on $\big(M^{\mathbb{N}_0}, \mathcal{B}(M)^{\otimes \mathbb{N}_0}\big)$ such that for all $B_0,\ldots,B_n \in \mathcal{B}(M)$,
\begin{align*}
&\mathbb{P}_{\lambda_0}\Big(
    \{ X = (X_i)_{i=0}^\infty \in M^{\mathbb{N}_0} \mid
    X_0 \in B_0, \ldots, X_n \in B_n \}
\Big)
\\ \nonumber
&= \int_{B_0} \lambda_0(dx_0) 
\int_{B_1} \mathcal{K}(x_0, dx_1) \cdots
\int_{B_n} \mathcal{K}(x_{n-1}, dx_n).
\end{align*}
Thus, a Markov chain can be identified with a random variable on $M^{\mathbb{N}_0}$ via coordinate projections.

Suppose that $\lambda$ is an invariant probability measure of the Markov chain $(X_n)_{n\geq 0}$.
The Markov chain is said to be stationary (or $\lambda$-stationary) if $\mathcal{L}(X_0) = \lambda$.
An invariant measure of a Markov kernel $\mathcal{K}$ is said to be ergodic if for every bounded and measurable function $f$ on $M$ the equation
\begin{equation*}
\int_{M}f(y)\mathcal{K}(x,dy)=f(x),
\end{equation*}
holds $\lambda$-a.s. for all $x\in M$, then $f$ is $\lambda$-a.s. a constant function.
We conclude this section by recalling the ergodic theorem for Markov chains on metric spaces.
\begin{theorem}[Ergodic Theorem for Markov Chains, \cite{benaim2022markov}]\label{markov ergodic thm}
Let $(X_n)_{n\geq 0}$ be a Markov chain on a Polish metric space $M$ with transition kernel $\mathcal{K}$, and let $\mathbb{P}_\lambda$ denote the probability measure on the path space $M^{\mathbb{N}_0}$ assuming the initial distribution is $X_0 \sim \lambda$.
If $\lambda$ is an ergodic probability measure for $\mathcal{K}$, then for any $h \in L^1(M, \lambda)$, we have
\begin{equation*}
\lim_{n\to\infty}\frac{1}{n}\sum_{k=0}^{n-1}h(X_k) = \int_M h(x)\,\lambda(dx) \quad \mathbb{P}_\lambda\text{-a.s.}
\end{equation*}
\end{theorem}

The ergodic theorem will play a central role in proving the strong consistency of the sample Fréchet mean in Section 4.

\section{Riemannian autoregression}

An $\mathbb{R}$-valued stochastic process $(X_n)_{n\geq 0}$ is said to be an autoregressive process of order one centred at $\mu\in \mathbb{R}$ if, for all $n \in \mathbb{N}$,
\begin{equation*}
X_n = (1-\phi)\mu + \phi X_{n-1} + \varepsilon_n,
\end{equation*}
where $\phi \in \mathbb{R}$ is a constant, and $(\varepsilon_n)_{n\geq 1}$ is a sequence of independent and identically distributed random variables.
This formulation makes sense in vector spaces, where addition and scalar multiplication are defined globally.
However, when extending this model to Riemannian manifolds, we must account for the absence of a global linear structure.
A natural approach is to utilise the tangent space at the Fréchet mean $\mu$, using the logarithmic map to lift points to the tangent space, and the exponential map to map them back.
Specifically, we realise a process as a Riemannian autoregression if its image under the log map is an ordinary autoregressive process.
Since the logarithmic and exponential maps are intrinsically defined through the geodesic structure induced by the Levi-Civita connection, this construction provides a geometric analogue of the ordinary autoregressive process.
To ensure that the log map is well-defined, i.e., single-valued, at each step, we constrain the update to lie within the global injectivity radius.

\subsection{Definition and basic properties}
We assume that $(\Omega,\mathcal{A},P)$ is a fixed probability space and let $M$ be a complete $d$-dimensional Riemannian manifold.
\begin{definition}[Riemannian autoregressive process]
Let $\mu\in M$ and $(\varepsilon_n)_{n\ge 1}$ be an independent and identically distributed sequence of non-degenerate $T_\mu M$-valued random variables satisfying $E[\varepsilon_n]=0$ and $E[\norm{\varepsilon_n}^2]<\infty$.
Let $X_0$ be an $M$-valued random variable, independent of $(\varepsilon_n)_{n\ge 1}$, and such that $E[d^2(X_0,\mu)]<\infty$.
An $M$-valued stochastic process $(X_n)_{n\ge 0}$ is called a Riemannian autoregressive process of order one, or simply Riemannian autoregressive process at $\mu$ with coefficient $\phi\in\mathbb R$, denoted $\mathrm{R\text{-}AR}(\mu,\phi)$, if there exists a radius $0<r\leq \injr{M}$ such that, with $B_0(r)$ denoting the open ball of radius $r$ centred at $0$ in $T_\mu M$,
\begin{enumerate}
\item $\log_{\mu}X_0 \in B_0(r)$, $P$-a.s.,
\item for all $n\ge 1$,
\begin{equation}\label{eq: R-AR def}
\log_{\mu}X_n = \phi \log_\mu(X_{n-1}) + \varepsilon_n\in B_{0}(r),
\qquad \text{$P$-a.s.}.
\end{equation}
\end{enumerate}
If $\injr{M}=\infty$, we adopt the conventions $B_\mu(\infty)=M$ and $B_0(\infty)=T_\mu M$.
\end{definition}
In principle, the model could be defined without second-moment assumptions by relying solely on the almost-sure containment of the process within a fixed geodesic ball. 
However, assuming $E[d^2(X_0, \mu)] < \infty$ and $E[\norm{\varepsilon_n}^2] < \infty$ simplifies the theory and, as we shall see in the following propositions, yields stronger results regarding its asymptotic stability.
Although the definition is presented here for a process of order one, with a scalar autoregressive parameter $\phi$, the formulation can be extended by considering higher-order dependence through additional lagged tangent vectors in $T_\mu M$, as well as by replacing $\phi$ with a linear operator on the tangent space.

We note that when $M=\mathbb{R}^d$, the proposed recursion corresponds to a restricted vector autoregressive model with coefficient matrix $\phi I_d$.
The recursion
\[
X_n=\exp_\mu\bigl(\phi\log_\mu X_{n-1}+\varepsilon_n\bigr)
\]
defines a local autoregressive model on any Riemannian manifold, since it is formulated entirely in the tangent space $T_\mu M$.
Whether this local model extends to a globally well-defined stochastic process depends on the geometry of the underlying manifold.
In particular, the existence of a finite radius $r<\infty$, such that $X_n\in B_\mu(r)$ $P$-a.s. for all $n\geq 0$
imposes strong restrictions on the parameters of the model.
If $M$ has finite diameter, in particular if $M$ is compact, then any Riemannian autoregressive process admitting such a finite containment radius must necessarily have uniformly bounded noise and satisfy $|\phi|<1$, reflecting the fact that a finite diameter implies a bounded injectivity radius.
In contrast, on Hadamard manifolds the injectivity radius is infinite, and the above local recursion extends globally (with $r=\infty$), in which case no such restriction on $\phi$ or on the noise is imposed.
We make this distinction precise in the following proposition.

\begin{proposition}\label{proposition: bounded injectivity radius and phi property}
Let $M$ be a Riemannian manifold. 
Suppose $(X_n)_{n\geq 0}$ is an $\mathrm{R\text{-}AR}(\mu,\phi)$ on $M$ and $X_0\in B_{\mu}(r_0)$ $P$-a.s.\ with $r_0<\infty$.
Then there exists $r<\infty$ such that $X_{n}\in B_{\mu}(r)$ $P$-a.s.\ for all $n\in\mathbb{N}$ if and only if $\abs{\phi}<1$ and there exists $\delta>0$ such that $\norm{\varepsilon_n}<\delta<\infty$ $P$-a.s.\ for all $n\in\mathbb{N}$.
\end{proposition}
\begin{proof}
Suppose that $\abs{\phi}<1$ and $\norm{\varepsilon_n}<\delta<\infty$ $P$-a.s.\ for all $n\in\mathbb{N}$. 
We want to show that there exists $r<\infty$ such that $X_n\in B_{\mu}(r)$ $P$-a.s.\ for all $n\in\mathbb{N}$.
Iteratively using the triangle inequality, we obtain
\begin{equation*}
\norm{\log_\mu X_n}\leq \abs{\phi}^n \norm{\log_{\mu}X_0} + \delta\sum_{k=0}^{n-1} \abs{\phi}^k
=
\abs{\phi}^n \norm{\log_{\mu}X_0} + \delta\frac{1-\abs{\phi}^n}{1-\abs{\phi}}.
\end{equation*}
If $\norm{\log_{\mu}X_0}\le r_0$ $P$-a.s., then for all $n\geq 1$,
\begin{equation*}
\norm{\log_{\mu}X_n}<\max\left\{r_0,\frac{\delta}{1-\abs{\phi}}\right\}
\quad P\text{-a.s.}
\end{equation*}
Consequently, $X_n\in B_\mu(r)$ $P$-a.s.\ for all $n\ge 1$, where $r=\max\left\{r_0,\frac{\delta}{1-\abs{\phi}}\right\}$.

Conversely, suppose that $X_n\in B_{\mu}(r)$ $P$-a.s., where $r<\infty$.
We want to show that $\abs{\phi}<1$ and there exists $\delta>0$ such that $\norm{\varepsilon_n}<\delta<\infty$ $P$-a.s.\ for all $n\in\mathbb{N}$.
Since $X_n$ is $\mathrm{R\text{-}AR}(\mu,\phi)$ and $X_n\in B_{\mu}(r)$, we have for all $n\in\mathbb{N}$
\begin{equation*}
\varepsilon_{n}=\log_{\mu}X_{n}-\phi\log_{\mu}X_{n-1}
\implies
\norm{\varepsilon_n}\leq R(1+\abs{\phi})<\infty \quad P\text{-a.s.},
\end{equation*}
where $R=\max\{r_0,r\}$.
So we can take $\delta=R(1+\abs{\phi})$.
Since $\varepsilon_n$ is non-degenerate in $T_\mu M$, there exists a unit vector
$u\in T_\mu M$ such that
\begin{equation*}
\var\left(\inpr{u}{\varepsilon_n}\right)=c>0,
\end{equation*}
for all $n\in\mathbb{N}$.
Then
\begin{equation*}
\abs{\inpr{u}{\log_{\mu}X_n}}\leq \norm{\log_{\mu}X_n}<R \quad P\text{-a.s.}
\end{equation*}
Consider the random variable $Y_n=\inpr{u}{\log_\mu X_n}$.
As $\abs{Y_n}<R$ $P$-a.s., we have $\var({Y_n})<R^2$.
By the i.i.d. assumption on the noise, $\var(\inpr{u}{\varepsilon_n})=c$, which is constant.
Furthermore, since $\varepsilon_n$ is independent of $X_{n-1}$ and has mean zero, the covariance between $Y_{n-1}$ and $\inpr{u}{\varepsilon_n}$ is exactly zero. Thus,
\begin{equation*}
\var({Y_n})
=
\phi^2 \var({Y_{n-1}}) + c
=
\phi^{2n}\var({Y_0}) + c\sum_{k=0}^{n-1}\abs{\phi}^{2k}
<
R^2.
\end{equation*}
Taking $n\to\infty$, we obtain
\begin{equation*}
\sum_{k=0}^{\infty}\abs{\phi}^{2k}<\infty,
\end{equation*}
and therefore $\abs{\phi}<1$.
\end{proof}

The second observation concerns the ergodicity of the process.
Specifically, we know that the Euclidean autoregressive process has a unique invariant law, if and only if $\abs{\phi}<1$.
For an R-AR process, we have the same situation.

\begin{proposition}\label{proposition: unique limiting distn and ergodicity of R-AR}
Let $(X_n)_{n\geq 0}$ be an R-AR$(\mu,\phi)$ process with noise
$(\varepsilon_n)_{n\geq 1}$, where each $\varepsilon_n$ has probability law
$\eta$, and let $r$ be as in the definition of the process. Assume that for all $x\in B_\mu(r)$
\begin{equation*}
\eta\set{\varepsilon\in T_\mu M}
{\phi\log_\mu x+\varepsilon\in B_0(r)}=1.
\end{equation*}
Then $(X_n)_{n\geq 0}$ defines a Markov chain on $B_\mu(r)$ with transition kernel
\begin{equation*}
\mathcal{K}(x,B)=\eta\set{\varepsilon\in T_{\mu}M}
{\exp_{\mu}(\phi\log_{\mu}x+\varepsilon)\in B},
\end{equation*}
for all $x\in B_\mu(r)$ and Borel sets $B\subseteq B_\mu(r)$.
Furthermore, its transition kernel has a unique invariant law if and only if
$\abs{\phi}<1$.
\end{proposition}

\begin{proof}
Let $(\varepsilon_n)_{n\geq 1}$ be an i.i.d. sequence with $\varepsilon_n \sim \eta$, and define $F_\varepsilon(x) = \exp_\mu\bigl(\phi \log_\mu x + \varepsilon\bigr)$. 
The process $(X_n)_{n\geq 0}$ follows the recursion $X_n = F_{\varepsilon_n}(X_{n-1})$ where $X_0$ is independent of the noise.
Since $(\varepsilon_n)_{n\geq 1}$ is independent of $X_{0},\ldots,X_{n-1}$, $(X_n)_{n\geq 0}$ is a Markov chain.
The associated transition kernel is given by
\begin{equation*}
\mathcal{K}(x,B) = \eta\set{\varepsilon \in T_\mu M}{\exp_\mu\big(\phi \log_\mu x + \varepsilon\big)\in B}.
\end{equation*}
For each $n\geq 0$, let $Y_n = \log_\mu X_n \in T_\mu M$. By the model recursion, the process $(Y_n)_{n\geq 0}$ satisfies 
\begin{equation}\label{proof of invariant law proposition: Y recursion}
Y_n = \phi Y_{n-1} + \varepsilon_n, \quad \text{for all } n \geq 1.
\end{equation} 
Suppose first that $\abs{\phi}<1$. 
Equation \eqref{proof of invariant law proposition: Y recursion} defines a classical autoregressive process in $T_\mu M$, which admits a unique invariant probability measure. 
Since $\exp_\mu$ is a diffeomorphism on $B_0(r)$, this law is pushed forward to a unique invariant probability measure $\lambda$ for $(X_n)_{n\geq 0}$.

Conversely, let $\lambda$ be an invariant probability measure for the transition kernel $\mathcal{K}$, and consider the stationary version of $(X_n)_{n\geq 0}$ with invariant law $\lambda$.
Then $(Y_n)_{n\geq 0}$ is also stationary as a measurable function of $(X_n)_{n\geq 0}$.
Because the noise is non-degenerate, there exists a unit vector $v \in T_\mu M$ such that $E[\inpr{\varepsilon_n}{v}^2] > 0$. 
Taking the inner product of \eqref{proof of invariant law proposition: Y recursion} with $v$ and using the stationarity of expectations, we have
\begin{equation*}
(1 - \phi^2) E[\inpr{Y_n}{v}^2] = E[\inpr{\varepsilon_n}{v}^2].
\end{equation*}
This forces $1 - \phi^2 > 0$, hence $\abs{\phi} < 1$.
\end{proof}

The parameter $\mu$ can be interpreted as the Fréchet mean, or the “centre,” of the process $\{X_n\}_{n\ge 0}$. 
Recall the distinction between two notions of mean. 
A Karcher mean is a local minimiser of the Fréchet function, whereas the Fréchet mean is its global minimiser. 
In the R-AR model setting, the recursion yields
\begin{equation*}
E[\log_\mu X_n] = \phi^n E[\log_\mu X_0],
\end{equation*}
since the noise is centred. 
Clearly, if $E[\log_\mu X_0]=0$, then $E[\log_\mu X_n]=0$ for all $n$.
This ensures that $\mu$ is a critical point of the Fréchet function for each $n$, and hence a Karcher mean provided $\mu$ is indeed a local minimiser. 
Whether $\mu$ is also the Fréchet mean depends heavily on the geometry of the space.
For instance, on Hadamard manifolds the two notions coincide because the Fréchet function is globally convex.
On more general manifolds, they coincide whenever the law of $X_n$ is supported in a geodesically strongly convex ball.
The following proposition gives a sufficient condition for $\mu$ to be the unique Fréchet mean of the invariant law of $X_n$, i.e., when the stationary law has $\mu$ as its Fréchet mean.
\begin{proposition}\label{proposition: mu interpretation as Fréchet mean of the stationary law}
Suppose $M$ is a complete Riemannian manifold whose curvature is bounded above by $K_U\in\mathbb{R}$.
Let $(X_n)_{n\ge 0}$ be a stationary R-AR$(\mu,\phi)$ process with invariant law $\lambda$, and assume the Fréchet function of $\lambda$ is finite on $M$.
Assume $X_n$ is supported in $B_{\mu}(r^*)$, where
\begin{equation}\label{inequality: mu interpretation convexity radius bound inequality}
r^* < \frac{1}{2}\min\left\{\text{inj}(M), \frac{\pi}{\sqrt{K_U}}\right\},
\end{equation}
with the convention that $\frac{1}{\sqrt{K_U}}=\infty$ if $K_U\leq 0$.
If $M$ is a Hadamard manifold, the upper bound in (\ref{inequality: mu interpretation convexity radius bound inequality}) evaluates to $\infty$, and we take $r^* = \infty$ ($B_\mu(\infty)=M$). 
Then $\mu$ is the unique Fréchet mean of $\lambda$.
\end{proposition}

\begin{proof}
By the stationarity assumption, taking expectation with respect to $\lambda$ for both sides in $\log_{\mu}X_n=\phi \log_{\mu}X_{n-1}+\varepsilon_n$ gives $E[\log_\mu X_n]=\phi E[\log_{\mu}X_n]$. 
By Proposition \ref{proposition: unique limiting distn and ergodicity of R-AR}, $|\phi|<1$. 
Hence, the only solution is $E[\log_\mu X_n]=0$ for all $n\geq 0$. 
Thus $E[\log_\mu X_n]=0$ and $\mu$ is a critical point of the Fréchet function of $\lambda$.
If $r^* < \infty$, then $X_n\in B_\mu(r^*)$, and by Theorem 2.1 in \cite{afsari2011riemannian} the Fréchet function admits a unique critical point in $B_\mu(r^*)$, and therefore, $\mu$ is the Fréchet mean of $\lambda$.
If $r^*=\infty$, i.e., the unbounded-support case, and $M$ is a Hadamard manifold, the convexity of distance function gives the uniqueness of the critical points of the Fréchet function; see \cite{sturm2003probability} or \cite{pennec2006intrinsic}.
Therefore, $\mu$ is the unique Fréchet mean of $\lambda$.
\end{proof}

On the other hand, the parameter $\phi$ governs the asymptotic behaviour of the process.
If $|\phi|<1$, then the tangent process $\log_\mu X_n$ is an ordinary Euclidean autoregressive process admitting a unique stationary distribution.
The pushforward of this stationary law under $\exp_\mu$ yields a stationary distribution for the Riemannian autoregressive process, whose Fréchet mean is $\mu$.
In this regime, the process exhibits mean--reverting behaviour toward $\mu$.
In contrast, if $\abs{\phi}\geq 1$, then the process is well-defined for all $n\geq 0$ only when the manifold $M$ has infinite injectivity radius, as shown in Proposition~\ref{proposition: bounded injectivity radius and phi property}.
We illustrate the resulting behaviour by considering the hyperbolic plane as an example, highlighting the distinct roles played by the parameters $\mu$ and $\phi$.
\begin{example}
The hyperbolic plane $\mathbb{H}^2$ is a complete 2-dimensional Riemannian manifold whose sectional curvature is $-1$.
There are various geometric realisations of hyperbolic space, known as models of $\mathbb{H}^2$, all of which are mutually isometric \cite{Rie-Lee}.
In this example, we consider the Poincaré disk model of $\mathbb{H}^2$.
Consider the unit open ball as an embedded submanifold of $\mathbb{R}^2$,
\begin{equation*}
\mathbb{B}^2 = \set{x\in\mathbb{R}^2}{\norm{x}<1}.
\end{equation*}
Let $(x^1,x^2)$ be the standard coordinates on $\mathbb{B}^2$ with corresponding coordinate vector fields $\pd{}{x^1}$ and $\pd{}{x^2}$.
The hyperbolic metric in these coordinates is the Riemannian metric tensor
\begin{equation*}
g_{\mathbb{H}^2}=\frac{4}{(1-\norm{x}^2)^2}
\left( dx^1 \otimes dx^1 + dx^2 \otimes dx^2 \right),
\end{equation*}
where $dx^1$ and $dx^2$ form the basis of the dual space of $T_x\mathbb{B}^2$ 
corresponding to $\pd{}{x^1}$ and $\pd{}{x^2}$.
Equivalently, for each $x\in\mathbb{B}^2$ and all $v,w\in T_x\mathbb{B}^2$,
\begin{equation*}
g_{\mathbb{H}^2,x}(v,w)=\frac{4}{\left(1-\norm{x}^2\right)^2}\left( v^1 w^1 + v^2 w^2 \right).
\end{equation*}
The Poincaré disk or ball model of the hyperbolic plane $\mathbb{H}^2$ is the Riemannian manifold $(\mathbb{B}^2, g_{\mathbb{H}})$.
As a Riemannian manifold, $\mathbb{H}^2$ is geodesically complete and $\injr{(\mathbb{H}^2)}=\infty$.
The geodesics in the Poincaré disk consist of either Euclidean straight lines, when the geodesic passes through the centre of the disk, or Euclidean circular arcs that intersect the boundary circle orthogonally.
For explicit formulas for the exponential and logarithm maps we refer to \cite{lou2020differentiating}.

We illustrate the roles of $\mu$ and $\phi$ for an R-AR process in the Poincaré disk.
Let $\mu=(\tanh(0.25),0)\in \mathbb{B}^2$.
We consider the process $X_n=\exp_{\mu}\left(\phi\log_{\mu}X_{n-1}+\varepsilon_n\right)$, where the noise is uniformly sampled from a disk of radius $\delta=1.5$ centred at $0$ in $T_\mu M$.
See Fig.\ref{fig:traj} and \ref{fig:rays} for point evolution of an R-AR in hyperbolic space $\mathbb{H}^2$ for different values of $\phi$.
In the case of $\phi=1.2$, the drift away from $\mu$ is evident.
By contrast, if $\phi=0.8$ we see the mean-reversion behaviour of $X_n$ towards $\mu$.
Because the hyperbolic plane has constant negative curvature ($K= -1$), the radius bound on $r^*$ in Proposition~\ref{proposition: mu interpretation as Fréchet mean of the stationary law} becomes infinite. 
Consequently, we are guaranteed that $\mu$ is the unique Fréchet mean of the stationary law of $X_n$.

\begin{figure}[H]
    \centering
    \begin{subfigure}[b]{0.48\textwidth}
        \centering
        \includegraphics[width=\textwidth]{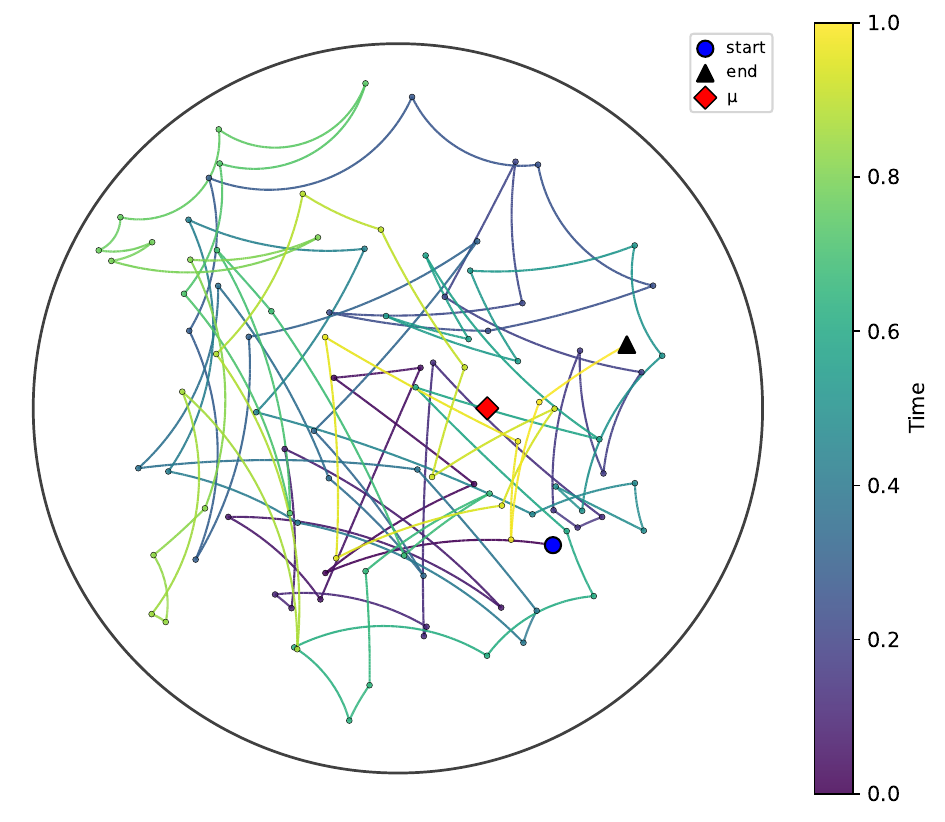}
        \caption{$\phi = 0.8$}
        \label{fig:traj_phi08}
    \end{subfigure}
    \hfill
    \begin{subfigure}[b]{0.48\textwidth}
        \centering
        \includegraphics[width=\textwidth]{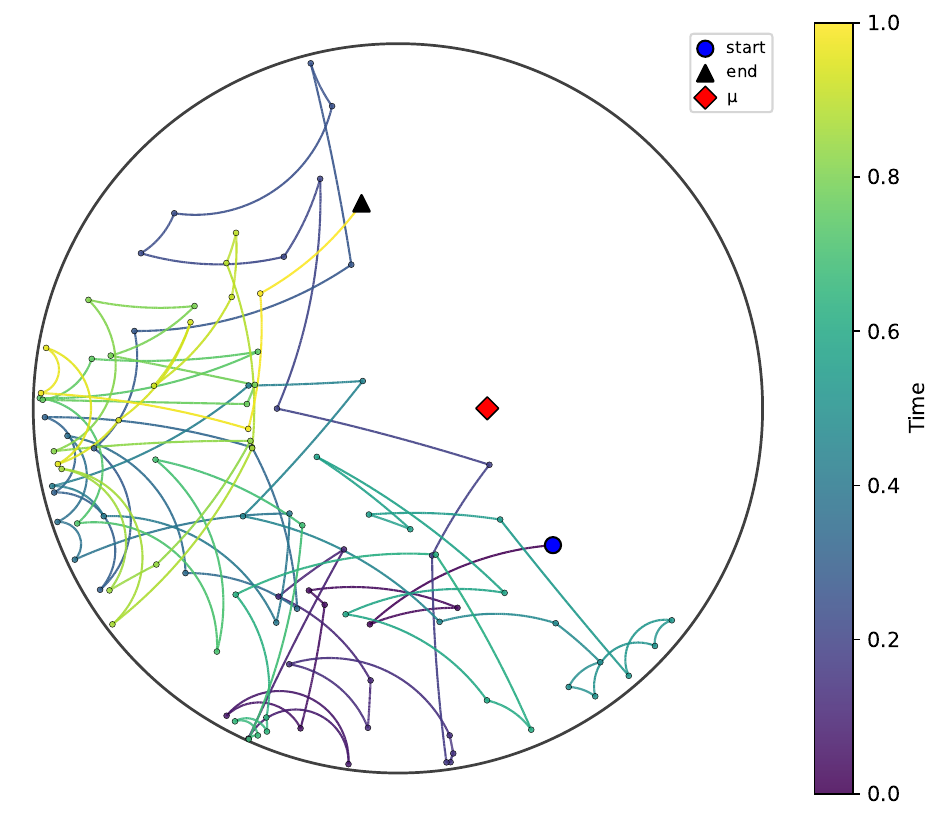}
        \caption{$\phi = 1.2$}
        \label{fig:traj_phi12}
    \end{subfigure}
    \caption{Trajectory evolution in the Poincaré disk model of $\mathbb{H}^2$ for two R-AR models with the same $\mu$ and different values of $\phi$.}
    \label{fig:traj}
\end{figure}

\begin{figure}[H]
    \centering
    \begin{subfigure}[b]{0.48\textwidth}
        \centering
        \includegraphics[width=\textwidth]{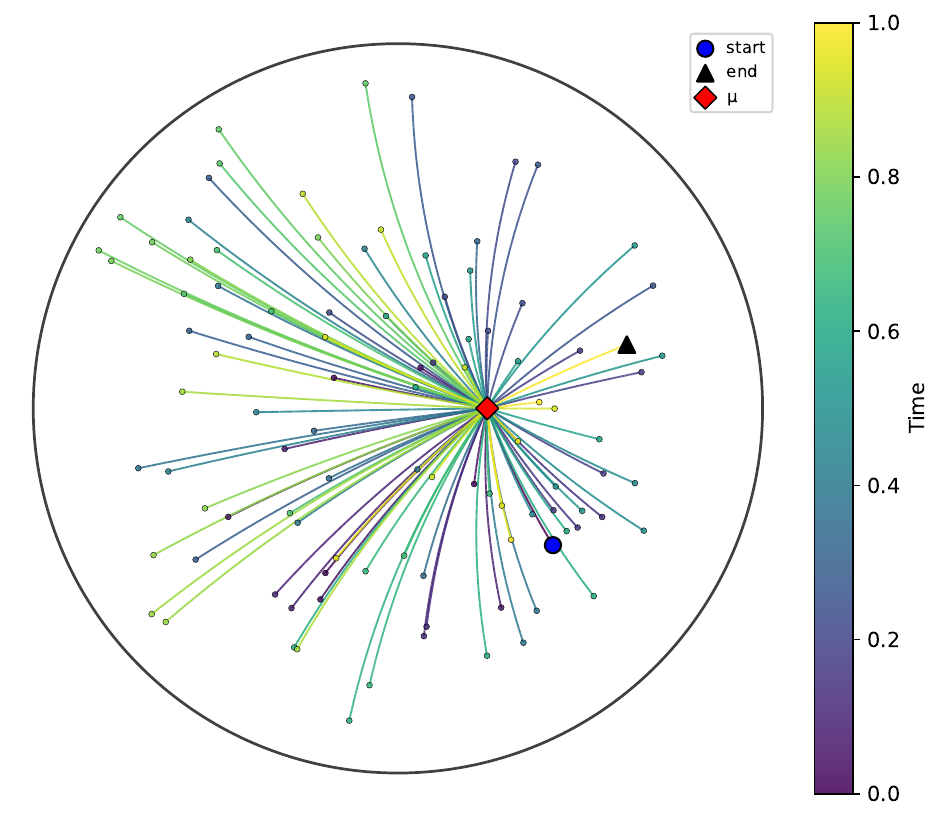}
        \caption{$\phi = 0.8$}
        \label{fig:rays_phi08}
    \end{subfigure}
    \hfill
    \begin{subfigure}[b]{0.48\textwidth}
        \centering
        \includegraphics[width=\textwidth]{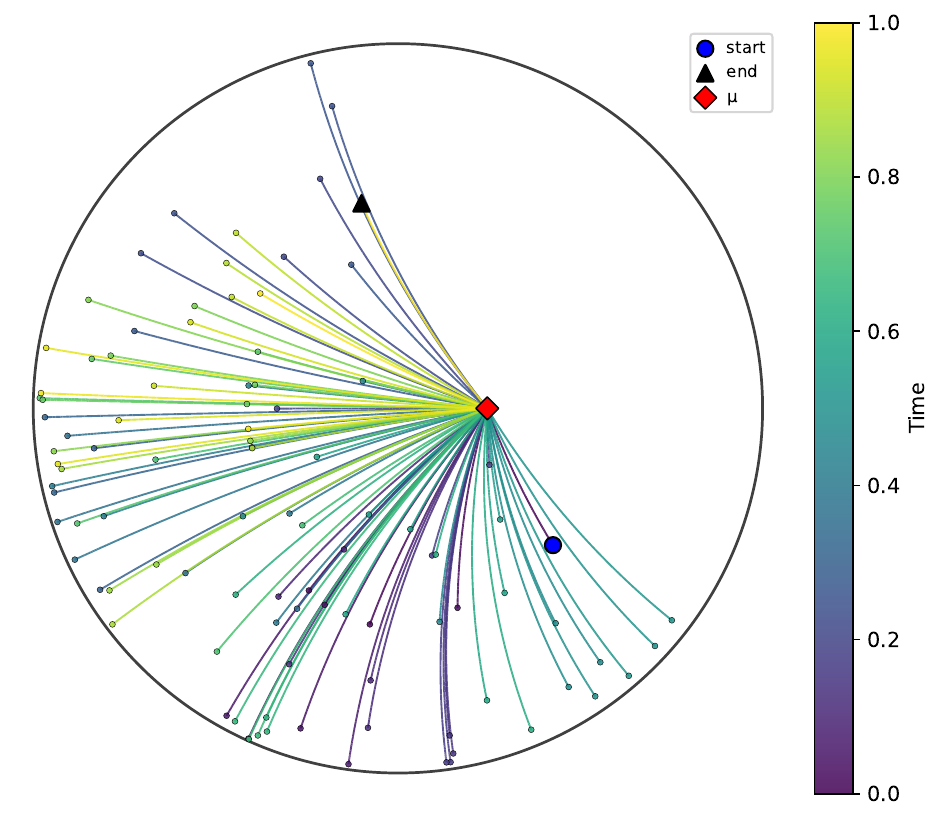}
        \caption{$\phi = 1.2$}
        \label{fig:rays_phi12}
    \end{subfigure}
    \caption{Geodesic rays from $\mu$ in the Poincaré disk model of $\mathbb{H}^2$ for two R-AR models with the same $\mu$ and different values of $\phi$.}
    \label{fig:rays}
\end{figure}
\end{example}

\subsection{Estimators}\label{sec: estimators subsection}
Consider a stationary R-AR$(\mu,\phi)$ process taking values in a Riemannian manifold $M$, and let $\{X_1,\ldots,X_n\}$ denote observed data points from the process. 
In this subsection, we propose estimators for the parameters $\mu$ and $\phi$ of the process.
A natural estimator for the parameter $\mu$ assuming it is the Fréchet mean of the process is the sample Fréchet mean, defined as
\begin{equation*}
\mu_n = \mathrm{arg}\min_{p \in M} \left( \frac{1}{n} \sum_{k=1}^n d^2(X_k, p) \right),
\end{equation*}
where $d(\cdot,\cdot)$ denotes the geodesic distance on the manifold $M$.
Next, we turn to the estimation of the autoregressive coefficient $\phi$. 
A direct approach to estimating $\phi$ is to solve for it from the model equation (\ref{eq: R-AR def}).
Under stationarity, the recursion in $T_\mu M$ gives
\begin{equation*}
\log_\mu X_n = \phi \log_\mu X_{n-1} + \varepsilon_n,
\end{equation*}
for all $n\geq 1$.
Taking inner products with $\log_\mu X_{n-1}$ and expectations yields
\begin{equation*}
E\left[\inpr{\log_\mu X_k}{\log_\mu X_{k-1}}\right]=\phi E\left[\norm{\log_\mu X_{k-1}}^2\right],
\end{equation*}
and therefore
\begin{equation*}
\phi=\frac{E\left[\inpr{\log_\mu X_k}{\log_\mu X_{k-1}}\right]}{E\left[\norm{\log_\mu X_{k-1}}^2\right]}.
\end{equation*}
Replacing expectations by sample averages leads to the estimator
\begin{equation*}
\phi_n=\frac{\sum_{k=2}^n\inpr{\log_{\mu_n} X_k}{\log_{\mu_n} X_{k-1}}}{\sum_{k=2}^n\norm{\log_{\mu_n} X_{k-1}}^2}.
\end{equation*}

The central question is the asymptotic behaviour of these estimators.
In the independence setting, strong consistency (almost sure convergence) of the sample Fréchet mean is established in \cite{bhattacharya2003large}.
However, the R-AR process generates dependent manifold-valued data, and these classical results cannot be applied directly. 
To address this, it is necessary to extend the theory to accommodate the dependence structure in the R-AR process. 
The next section develops a general convergence framework for sample Fréchet means under a dependence setting in metric spaces.
Using this framework, we establish strong consistency of both $\mu_n$ and $\phi_n$, with the R-AR process appearing as a special case of the general result.

\section{Asymptotic properties}
In this section, we prove the strong consistency of the R--AR parameter estimators. 
Our approach establishes a general convergence theorem in metric spaces, which we then apply to show the strong consistency of the estimator $\mu_n$ in Subsection \ref{sec: estimators subsection}.
The consistency of the autoregressive coefficient estimator $\phi_n$ is then established using additional geometric arguments.
In particular, we prove the strong consistency of the sample Fréchet mean for ergodic Markov chains in metric spaces, which is then used in showing the convergence of the proposed parameter estimators.

\subsection{Strong consistency of sample Fréchet means for Markov chains in metric spaces}
The objective of this subsection is to establish the strong consistency of the sample Fréchet mean set for a Markov chain with values in a proper metric space.
The strong law of large numbers for Fréchet means that we prove below can also be deduced, in a more general setting, from the results of \cite{jaffe2024fr}.
Nevertheless, we include a direct proof following the argument of \cite{bhattacharya2003large}, with the ergodic theorem replacing the strong law of large numbers used under independence, since this proof strategy is also used in establishing the convergence of our parameter estimators.
To achieve this, we first present a preliminary lemma whose conclusions will be used in the proofs of convergence theorems.
The claims of this lemma are proved in \cite{bhattacharya2003large} for manifolds, but the same arguments extend directly to general metric spaces. 
For the main theorem, we assume that the metric space is proper. A metric space is said to be proper if its closed balls are compact. 
Equivalently, proper metric spaces have the Heine-Borel property: closed and bounded subsets of $M$ are compact.
In particular, every proper metric space is complete and separable.
Therefore, in the setting of the main theorem, the Markov-chain ergodic theorem applies.

\begin{lemma}\label{lemma: properties of Fréchet function}
Let $(M,d)$ be a metric space.
Suppose that $X:\Omega \to M$ is a random variable with probability law $\lambda$.
Suppose that the Fréchet function $\mathcal{F}$ of $X$ is finite at some point $p_0 \in M$.
Then
\begin{enumerate}
\item $\mathcal{F}$ is finite on $M$,
\item $\mathcal{F}$ is continuous on $M$,
\item if $M$ is proper, then the Fréchet mean set of $X$ is a nonempty compact subset of $M$.
\end{enumerate}
\end{lemma}

\begin{proof}
See Appendix \ref{sec:appendix:proofs}
\end{proof}
For the next result, we impose the properness assumption on $M$. 
Recall that a metric space is proper if its closed balls are compact.
In particular, proper metric spaces are complete and separable.
Thus, under this assumption, the Markov-chain ergodic theorem applies; see \cite{benaim2022markov}.
\begin{theorem}[Almost sure uniform convergence on compact sets]
\label{thm: a.s. uniform convergence Fréchet functions on compact sets}
Let $(M,d)$ be a separable metric space and let
$(X_n)_{n\ge0}$ be an $M$-valued Markov chain defined on a probability space $(\Omega,\mathcal{A}, P)$ with an ergodic probability measure $\lambda$.
Assume that $X_0\sim\lambda$ and that the Fréchet function
\begin{equation*}
\mathcal F(p_0)=\int_M d^2(x,p_0)\,\lambda(dx)
\end{equation*}
is finite for some $p_0\in M$.
Define
\begin{equation*}
\mathcal F_n(p)=\frac{1}{n+1}\sum_{k=0}^{n} d^2(X_k,p).
\end{equation*}
Then, for every compact set $S\subset M$,
\begin{equation*}
\sup_{p\in S}\abs{\mathcal F_n(p)-\mathcal F(p)}
\xrightarrow[n\to\infty]{} 0
\qquad  P\text{-a.s.}
\end{equation*}
\end{theorem}

\begin{proof}
Let $\Omega' = M^{\mathbb{N}_0}$ be the trajectory space equipped with the product Borel $\sigma$-algebra. Let $X: \Omega \to \Omega'$ be the mapping defined by $X = (X_0, X_1, \dots)$, and let $\mathbb{P}_\lambda = P \circ X^{-1}$ be the pushforward probability measure on $\Omega'$. We first establish the convergence for trajectories $\omega = (\omega_0, \omega_1, \dots) \in \Omega'$ under the measure $\mathbb{P}_\lambda$.

Fix a compact set $S\subset M$.
From Lemma~\ref{lemma: properties of Fréchet function},
$\mathcal F$ is finite on $M$.
Note that for all $p,q\in S$
\begin{equation*}
\abs{\mathcal{F}_n(p, \omega)-\mathcal{F}(p)}
\le
\abs{\mathcal{F}_n(p, \omega)-\mathcal{F}_n(q, \omega)}
+
\abs{\mathcal{F}_n(q, \omega)-\mathcal{F}(q)}
+
\abs{\mathcal{F}(q)-\mathcal{F}(p)}.
\end{equation*}
We estimate the three terms on the right-hand side separately and combine these estimates to obtain a bound uniform in $p\in S$.
For each $p\in S$, define
\begin{equation*}
f_p(\omega)=d(\omega_0,p),
\end{equation*}
where $\omega=(\omega_0,\omega_1,\dots)\in\Omega'$.
Since $\mathcal F(p)<\infty$, it follows that
$f_p\in L^1(\Omega',\mathbb P_\lambda)$.
By Theorem~\ref{markov ergodic thm}, for each fixed $p\in S$
\begin{equation*}
\frac{1}{n+1}\sum_{k=0}^{n}d(\omega_k,p)
\xrightarrow[n\to\infty]{\mathbb P_\lambda\text{-a.s.}}
E_\lambda[d(X_0,p)].
\end{equation*}
Fix $s_0\in S$.
Since $S$ is compact, $\operatorname{diam}(S)<\infty$.
For $p,q\in S$,
\begin{equation*}
\abs{d^2(x,p)-d^2(x,q)}=\abs{d(x,p)-d(x,q)}(d(x,p)+d(x,q))\leq d(p,q)(d(x,p)+d(x,q)).
\end{equation*}
Using the triangle inequality,
\begin{equation*}
d(x,p)
\leq d(x,s_0)+d(p,s_0)\leq d(x,s_0)+\operatorname{diam}(S),
\end{equation*}
and similarly for $d(x,q)$. 
Hence
\begin{equation*}
\abs{d^2(x,p)-d^2(x,q)}
\leq
2\,d(p,q)\left(d(x,s_0)+\operatorname{diam}(S)\right).
\end{equation*}
Therefore,
\begin{equation*}
\abs{\mathcal F_n(p,\omega)-\mathcal F_n(q,\omega)}
\leq
2d(p,q)
\left(
\frac{1}{n+1}\sum_{k=0}^{n}d(\omega_k,s_0)
+
\operatorname{diam}(S)
\right).
\end{equation*}

By the ergodic theorem applied to $d(\omega_0,s_0)$, there exists a set $A \subset \Omega'$ such that $\mathbb{P}_{\lambda}(A)=1$ and for each $\omega\in A$ there is $N_0(\omega)$ with
\begin{equation*}
\frac{1}{n+1}\sum_{k=0}^{n}d(\omega_k,s_0)
\leq
E_\lambda[d(X_0,s_0)] + 1
\quad
\text{for all } n\ge N_0(\omega).
\end{equation*}

Hence for $n\ge N_0(\omega)$,
\begin{equation*}
\abs{\mathcal F_n(p,\omega)-\mathcal F_n(q,\omega)}
\leq
\alpha\, d(p,q),
\end{equation*}
for all $p,q\in S$, where
\begin{equation*}
\alpha=
2\big(E_\lambda[d(X_0,s_0)] + 1
+
\operatorname{diam}(S)\big).
\end{equation*}

Next, for each fixed $p\in S$, using the ergodic theorem for $d^2(\omega_0,p)$,
\begin{equation}\label{proof: pointwise convergence of sample Fréchet function}
\mathcal F_n(p,\omega)\to\mathcal F(p)
\quad
\mathbb P_\lambda\text{-a.s.}
\end{equation}

Fix $\varepsilon>0$.
Since $S$ is compact, $\mathcal{F}$ is uniformly continuous on $S$, and thus, there exists $\delta_1>0$ such that
\begin{equation*}
d(p,q)<\delta_1
\implies
\abs{\mathcal F(p)-\mathcal F(q)}<\frac{\varepsilon}{3}.
\end{equation*}
Let
$\delta=\min\left\{\delta_1,\frac{\varepsilon}{3\alpha}\right\}$.
Since $S$ is compact, there exist finitely many points $q_1,\dots,q_m\in S$ such that $S\subset \bigcup_{j=1}^{m} B(q_j,\delta)$.
For each $j$, define
\begin{equation*}
B_{q_j}=\set{\omega\in\Omega'}{\mathcal F_n(q_j,\omega)\to\mathcal F(q_j)}.
\end{equation*}
Then, by \eqref{proof: pointwise convergence of sample Fréchet function}, $\mathbb P_\lambda(B_{q_j})=1$, for all $j=1,\ldots,m$.
Define
\begin{equation*}
C=A\cap\bigcap_{j=1}^{m}B_{q_j}.
\end{equation*}
Then $\mathbb P_\lambda(C)=1$.
Fix $\omega\in C$.
For each $j$, there exists $N_j(\omega)$ such that
\begin{equation*}
n\ge N_j(\omega)
\implies
|\mathcal F_n(q_j,\omega)-\mathcal F(q_j)|<\frac{\varepsilon}{3}.
\end{equation*}
Let $N(\omega)=\max_{0\leq i \leq m}\{N_i(\omega)\}$.
For $n\ge N(\omega)$ and any $p\in S$, choose $q_j$ such that $d(p,q_j)<\delta$. Then
\begin{equation*}
\begin{aligned}
|\mathcal F_n(p,\omega)-\mathcal F(p)|
&\le
|\mathcal F_n(p,\omega)-\mathcal F_n(q_j,\omega)|
+
|\mathcal F_n(q_j,\omega)-\mathcal F(q_j)|
+
|\mathcal F(q_j)-\mathcal F(p)| \\
&<
\frac{\varepsilon}{3}
+
\frac{\varepsilon}{3}
+
\frac{\varepsilon}{3}
=
\varepsilon.
\end{aligned}
\end{equation*}

Therefore, the set of trajectories 
\begin{equation*}
E = \left\{ \omega \in \Omega' : \sup_{p\in S} |\mathcal F_n(p,\omega)-\mathcal F(p)| \to 0 \right\}
\end{equation*}
has measure $\mathbb P_\lambda(E) = 1$. 
By the definition of the pushforward measure, we have
\begin{equation*}
P(X^{-1}(E)) = \mathbb{P}_\lambda(E) = 1.
\end{equation*}
Since $X^{-1}(E) = \left\{ \tilde{\omega} \in \Omega : \sup_{p\in S} |\mathcal F_n(p, \tilde{\omega}) - \mathcal F(p)| \to 0 \right\}$, the uniform convergence holds $P$-a.s.
\end{proof}

The next theorem is the main result of this section and is stated under the assumption that $(M,d)$ is a proper metric space, that is, a metric space in which closed balls are compact.
The result concerns the estimation of the Fréchet mean set.
Since this is a set-valued object, convergence is understood in terms of the distance from points in the sample Fréchet mean set to the Fréchet mean set.
For a subset $U \subset M$ and $x \in M$, define
\begin{equation*}
d(x,U) = \inf_{y \in U} d(x,y).
\end{equation*}
For a sequence of subsets $(U_n)_{n\geq 1}$ of $M$, the convergence result below is expressed through the quantity
\begin{equation*}
\sup_{x \in U_n} d(x,U),
\end{equation*}
which measures how far the elements of $U_n$ may be from the set $U$.
In the literature, this is referred to as the excess of $U_n$ over $U$, and it is closely related to the Hausdorff metric; see \cite{barbati1994hausdorff} or \cite{evans2024limit}.

\begin{theorem}[Strong consistency of sample Fréchet mean set]\label{thm: strong consistency of Fréchet mean set for markov chains}
Let $(M,d)$ be a proper metric space.
Let $(X_n)_{n\ge0}$ be an $M$-valued Markov chain defined on a probability space $(\Omega, \mathcal{A}, P)$ with an ergodic probability measure $\lambda$, and assume that $ X_0\sim\lambda $.
Assume that the Fréchet function
\begin{equation*}
\mathcal F(p_0)=\int_M d^2(x,p_0)\lambda(dx),
\end{equation*}
is finite for some $p_0\in M$.
Let $\mathcal U$ denote the set of minimisers of $\mathcal F$, and $\mathcal U_n$ the set of minimisers of the sample Fréchet function
\begin{equation*}
\mathcal F_n(p)=\frac1n\sum_{k=0}^{n-1} d^2(X_k,p).
\end{equation*}
Then
\begin{equation*}
\sup_{q\in\mathcal U_n}d(q,\mathcal{U})\xrightarrow[n\to\infty]{P\text{-a.s.}}0.
\end{equation*}
\end{theorem}
\begin{proof}
We show that for all $\varepsilon>0$, there exists $N\in\mathbb{N}$ such that for all $n>N$, we have
\begin{equation*}
d(q,\mathcal{U})<\varepsilon,
\end{equation*}
for all $q\in\mathcal{U}_{n}$.
First note that by Lemma \ref{lemma: properties of Fréchet function}, we have $\mathcal{U}\neq \emptyset$, as the space is proper.
Fix $p^*\in\mathcal{U}$ and let
\begin{equation*}
m=\min_{p\in M}\mathcal F(p)=\mathcal F(p^*).
\end{equation*}

Fix $\varepsilon>0$ and define the $\varepsilon$-neighbourhood of $\mathcal U$ by
\begin{equation*}
\mathcal U_\varepsilon=\{p\in M:\ d(p,\mathcal U)<\varepsilon\}.
\end{equation*}

If $M$ is bounded, then it is compact, since it is proper.
Furthermore, $M\setminus \mathcal{U}_\varepsilon$ is compact and disjoint from $\mathcal{U}$.
By continuity of $\mathcal{F}$, we may define
\begin{equation*}
\alpha=\min_{p\in M\setminus \mathcal{U}_\varepsilon}
\mathcal{F}(p) - m
>0.
\end{equation*}

By Theorem~\ref{thm: a.s. uniform convergence Fréchet functions on compact sets}, there exists $N_1\in\mathbb{N}$ such that for all $n>N_1$,
\begin{equation*}
\sup_{p\in M}\abs{\mathcal{F}_n(p)-\mathcal{F}(p)}<\alpha/2.
\end{equation*}

Let $q\in\mathcal{U}_n$.
Then for all $p\in\mathcal{U}$,
\begin{equation*}
\mathcal{F}_n(q)\le \mathcal{F}_n(p)
< \mathcal{F}(p)+\alpha/2
= m+\alpha/2.
\end{equation*}
Hence,
\begin{equation*}
\mathcal{F}(q)
< \mathcal{F}_n(q)+\alpha/2
< m+\alpha.
\end{equation*}
By definition of $\alpha$, this implies $q\notin M\setminus \mathcal{U}_\varepsilon$.
Therefore, for all $n>N=N_{1}$ and $q\in\mathcal{U}_n$, $d(q,\mathcal{U})<\varepsilon$.
Suppose that $M$ is unbounded.
Using the reverse triangle inequality,
\begin{equation*}
\abs{d(p,p^*)-d(p^*,x)}\leq d(x,p),
\end{equation*}
for all $x,p \in M$.
Squaring and integrating with respect to $\lambda$ yields
\begin{equation*}
d^2(p,p^*)+\mathcal F(p^*)-2d(p,p^*)\mathcal F_{1}(p^*)\leq \mathcal F(p),
\end{equation*}
where $\mathcal F_{1}(p^*)=\int_M d(x,p^*)\lambda(dx)<\infty$.
Since $\mathcal F(p^*)=m$, we have
\begin{equation*}
d^2(p,p^*)-2d(p,p^*)\mathcal F_{1}(p^*)\leq \mathcal F(p)-m.
\end{equation*}

Thus $\mathcal F(p)>m$ whenever $d(p,p^*)>2\mathcal F_{1}(p^*)$.
Fix $\eta>\max\{\varepsilon,2\mathcal F_{1}(p^*)\}$ and define the closed $\eta$-neighbourhood of $\mathcal U$ by
\begin{equation*}
\mathcal U_\eta=\{p\in M:\ d(p,\mathcal U)\le\eta\}.
\end{equation*}

Then for all $p\notin\mathcal U_\eta$ we have $d(p,\mathcal U)>\eta$, and hence $d(p,p^*)>\eta$ since $p^*\in\mathcal U$. Therefore,
\begin{equation*}
\mathcal F(p)-m
>
\eta^2-2\eta\mathcal F_1(p^*)
=:\alpha_1>0.
\end{equation*}
Hence $\mathcal F(p)>m+\alpha_1$ for all
$p\notin\mathcal U_\eta$.

Consider the set
\begin{equation*}
\mathcal{U}_{\varepsilon,\eta}
=\{p\in M:\ \varepsilon\le d(p,\mathcal U)\le\eta\}.
\end{equation*}
Since $M$ is proper and $\mathcal U$ is compact, $\mathcal{U}_{\varepsilon,\eta}$ is compact.
By continuity of $\mathcal{F}$ we have
\begin{equation*}
\alpha_2=\min_{p\in \mathcal{U}_{\varepsilon,\eta}}\mathcal{F}(p)-m>0.
\end{equation*}

Set $\alpha:=\min\{\alpha_1,\alpha_2\}>0$.
Applying Theorem~\ref{thm: a.s. uniform convergence Fréchet functions on compact sets} on the compact set $\mathcal{U}_{\eta}$, there exists $N_2\in\mathbb{N}$ such that for all $n>N_2$,
\begin{equation*}\label{ineq: sup over u_eta < alpha/2}
\sup_{p\in \mathcal{U}_{\eta}}\abs{\mathcal{F}_n(p)-\mathcal{F}(p)}<\alpha/2.
\end{equation*}

Let $q\in\mathcal{U}_n$. 
By the reverse triangle inequality, $d^2(x,q) \ge d^2(q,p^*) - 2d(q,p^*)d(x,p^*) + d^2(x,p^*)$ for any $x \in M$. 
So
\begin{equation*}
\mathcal{F}_n(q) - \mathcal{F}_n(p^*) \ge d(q,p^*) \left( d(q,p^*) - \frac{2}{n}\sum_{k=0}^{n-1} d(X_k,p^*) \right).
\end{equation*}
By the ergodic theorem, $\frac{1}{n}\sum_{k=0}^{n-1} d(X_k,p^*) \to \mathcal{F}_1(p^*)$ almost surely. 
Since we chose $\eta > 2\mathcal{F}_1(p^*)$, there exists $N_3 \in \mathbb{N}$ such that for all $n > N_3$, we have $\frac{2}{n}\sum_{k=0}^{n-1} d(X_k,p^*) < \eta$. 

If $q \notin \mathcal{U}_\eta$, then $d(q,p^*) \ge \eta$. 
Therefore, for $n > N_3$,
\begin{equation*}
\mathcal{F}_n(q) - \mathcal{F}_n(p^*) \ge d(q,p^*) \left( d(q,p^*) - \frac{2}{n}\sum_{k=0}^{n-1} d(X_k,p^*) \right) > d(q,p^*) (d(q,p^*) - \eta) \ge 0.
\end{equation*}
Thus $\mathcal{F}_n(q) > \mathcal{F}_n(p^*)$, which contradicts the minimality of $q$. Hence $q\in\mathcal U_\eta$ for all $n > N_3$.

Since $p^*\in\mathcal{U}\subset \mathcal{U}_\eta$, we have
\begin{equation*}
\mathcal F_n(q)\le \mathcal F_n(p^*)
< \mathcal F(p^*)+\alpha/2
= m+\alpha/2.
\end{equation*}
Hence,
\begin{equation*}
\mathcal{F}(q)
< \mathcal{F}_n(q)+\alpha/2
< m+\alpha.
\end{equation*}
By definition of $\alpha$, this implies $q\notin \mathcal{U}_{\varepsilon,\eta}$.
Therefore, $q\in\mathcal{U}_{\varepsilon}$, that is $d(q,\mathcal{U})<\varepsilon$.
So for all $n>N=\max\{N_2,N_3\}$ we have $d(q,\mathcal{U})<\varepsilon$, for all $q\in\mathcal{U}_{n}$.
\end{proof}

\subsection{Convergence of R-AR estimators}
We now apply Theorem \ref{thm: strong consistency of Fréchet mean set for markov chains} to the Riemannian autoregressive model to estimate $\mu$ in the stationary setup.
Under the condition $|\phi|<1$, the process admits a unique invariant ergodic distribution, and hence satisfies the stochastic assumptions of the strong consistency theorem.
\begin{theorem}
Let $M$ be a complete Riemannian manifold and let $(X_n)$ be a stationary
R-AR$(\mu,\phi)$ process with $\abs{\phi}<1$.
Let $\lambda$ denote its stationary probability measure and assume that
$\mu$ is the unique Fréchet mean of $\lambda$.
Then $\mu_n \xrightarrow{P\text{-a.s.}} \mu$ as $n\to\infty$.
\end{theorem}

\begin{proof}
By the Hopf-Rinow theorem, completeness of $M$ implies $M$ is a proper metric space.
Proposition~\ref{proposition: unique limiting distn and ergodicity of R-AR} establishes that $(X_n)$ is a Markov chain and, because $\abs{\phi}<1$, guarantees it admits a unique invariant probability measure $\lambda$.
Since unique invariant measures of Markov kernels are also ergodic, Proposition 4.30 \cite{benaim2022markov}, $\lambda$ is ergodic.
Thus the assumptions of
Theorem~\ref{thm: strong consistency of Fréchet mean set for markov chains}
are satisfied, and the sample Fréchet mean set converges $P$-a.s. to the Fréchet mean set of $\lambda$.
Since $\mu$ is the unique Fréchet mean of $\lambda$, the limit set reduces to the singleton $\{\mu\}$,
and therefore $\mu_n \xrightarrow{P\text{-a.s.}} \mu$.
\end{proof}
We now address estimation of the autoregressive parameter $\phi$ in the Riemannian autoregressive model.
In Subsection \ref{sec: estimators subsection}, we proposed the estimator
\begin{equation*}
\phi_n
=
\frac{
\sum_{k=1}^n\inpr{\log_{\mu_n} X_k}{\log_{\mu_n} X_{k-1}}
}{
\sum_{k=1}^n
\norm{\log_{\mu_n} X_{k-1}}^2
}.
\end{equation*}
Before we state and prove the strong consistency of $\phi_n$, we need the following lemma.
\begin{lemma}\label{lemma: bound on the inner product of logs}
Let $M$ be a geodesically complete Riemannian manifold, and fix a base point $\mu \in M$.  
Suppose the sectional curvatures of $M$ are bounded above by $K_U$ and below by $K_L$, for $K_U,K_L\in\mathbb{R}$. 
Let $B_{\mu}(r^*)$ be a geodesic ball of radius
\begin{equation*}
r^* < \frac{1}{2}\min\left\{\operatorname{inj}(M), \frac{\pi}{\sqrt{K_U}}\right\},
\end{equation*}
where $\frac{1}{\sqrt{K_U}}$ is interpreted as $\infty$ if $K_U \leq 0$. 
Let $S \subset B_\mu(r^*)$ be a compact subset.  
Then there exists $C_S > 0$ such that for all $p,q\in S$ and $x,y\in B_\mu(r^*)$,
\begin{equation*}
\abs{\inpr{\log_p x}{\log_p y} - \inpr{\log_q x}{\log_q y}}
\le
C_S d(p,q)\left(1 + d^2(\mu,x) + d^2(\mu,y)\right).
\end{equation*}
\end{lemma}
\begin{proof}
See Appendix \ref{sec:appendix:proofs}.
\end{proof}
\begin{theorem}
Suppose $M$ is a complete Riemannian manifold whose curvature is bounded above by $K_{U}$ and below by $K_L$, $K_U,K_L\in\mathbb{R}$.
Let $(X_n)$ be a stationary R-AR$(\mu,\phi)$ process with $|\phi|<1$,
and let $\lambda$ denote its unique invariant probability measure.
Assume that $\lambda$ is supported in a geodesic ball $B_{\mu}(r^*)$, where
\begin{equation*}
r^*<\frac{1}{2}\min\left\{\injr{M}, \frac{\pi}{\sqrt{K_U}}\right\},
\end{equation*}
(adopting the convention that $\frac{1}{\sqrt{K_U}}=\infty$ if $K_U \le 0$).
If $M$ is a Hadamard manifold, we take $r^* = \infty$ with $B_{\mu}(\infty)=M$.
Assume the Fréchet function of $\lambda$ is finite at some point.
Then 
\begin{equation*}
\phi_n
=
\frac{
\sum_{k=1}^n
\inpr{\log_{\mu_n}X_k}{\log_{\mu_n}X_{k-1}}
}{
\sum_{k=1}^n
\norm{\log_{\mu_n}X_{k-1}}^2
}
\xrightarrow[n\to\infty]{P\text{-a.s.}}
\phi .
\end{equation*}
\end{theorem}

\begin{proof}
From Proposition \ref{proposition: mu interpretation as Fréchet mean of the stationary law}, $\mu$ is the unique Fréchet mean of $\lambda$.
For any $n\in\mathbb{N}$, the point $\mu_n$ lies in $B_{\mu}(r^*)$, this follows by Theorem 2.1 in \cite{afsari2011riemannian}.
Since the geodesic ball $B_{\mu}(r^*)$ is strongly convex, for all $X_1,\ldots,X_n$, the expression $\log_{\mu_n}X_{k}$ is well-defined for all $k=1,\ldots,n$.

By strong consistency of $\mu_n$, we have $\mu_n \to \mu$ $P$-a.s. as $n\to\infty$.
Thus, there exists a compact set $S \subset B_\mu(r^*)$ containing $\mu$ such that $\mu_n \in S$ for all sufficiently large $n$.

For the denominator, note that $\frac{1}{n}\sum_{k=0}^{n-1}\norm{\log_{p}X_{k}}^2 = \mathcal{F}_n(p)$
and the triangle inequality gives
\begin{equation*}
\abs{\mathcal{F}_n(\mu_n) - \mathcal{F}(\mu)} 
\le 
\sup_{p \in S} \abs{\mathcal{F}_n(p) - \mathcal{F}(p)} 
+ 
\abs{\mathcal{F}(\mu_n) - \mathcal{F}(\mu)}.
\end{equation*}
By Theorem~\ref{thm: a.s. uniform convergence Fréchet functions on compact sets}, the uniform convergence on compact sets yields $\sup_{p \in S} \abs{\mathcal{F}_n(p) - \mathcal{F}(p)} \to 0$ $P$-a.s. as $n\to\infty$.
By the continuity of $\mathcal{F}$, we have $\abs{\mathcal{F}(\mu_n) - \mathcal{F}(\mu)} \to 0$. 
Therefore, the denominator converges $P$-a.s. to $\mathcal{F}(\mu) = E\left[\norm{\log_{\mu}X_0}^2\right]$.

For the numerator, define the sequence of functions
\begin{equation*}
H_n(p) = \frac{1}{n}\sum_{k=1}^n \inpr{\log_p X_k}{\log_p X_{k-1}}
\end{equation*}
and let $H(p) = E\left[\inpr{\log_p X_1}{\log_p X_0}\right]$. 
Lemma \ref{lemma: bound on the inner product of logs} gives
\begin{equation}\label{proof: estimate for inner product in the proof a.s. convergence of phi_n}
\abs{\inpr{\log_p x}{\log_p y} - \inpr{\log_q x}{\log_q y}} 
\le 
C_S d(p,q) \left(1 + d^2(\mu, x) + d^2(\mu, y)\right),
\end{equation}
where $C_S$ is a finite constant depending on $S$.
Let $G(x,y)= 1 + d^2(\mu, x) + d^2(\mu, y)$.
Because the Fréchet function is finite at $\mu$, we have $E[d^2(\mu, X_0)] < \infty$, which ensures the right-hand side of (\ref{proof: estimate for inner product in the proof a.s. convergence of phi_n}) is integrable with respect to $\lambda$. 
By the triangle inequality, for all $p,q \in S$,
\begin{equation*}
\abs{H_n(p) - H_n(q)} 
\le 
C_Sd(p,q)\frac{1}{n}\sum_{k=1}^n G(X_k, X_{k-1}).
\end{equation*}
Using exactly the same $\varepsilon/3$ argument in the proof of Theorem \ref{thm: a.s. uniform convergence Fréchet functions on compact sets} we obtain
\begin{equation*}
\sup_{p\in S}
\abs{H_n(p)-H(p)}
\xrightarrow[n\to\infty]{P\text{-a.s.}}0.
\end{equation*}
Furthermore, the triangle inequality and the almost sure convergence of $\mu_n$ to $\mu$,
\begin{equation*}
\abs{H_n(\mu_n) - H(\mu)}
\leq \sup_{p \in S} \abs{H_n(p) - H(p)} + \abs{H(\mu_n) - H(\mu)} 
\xrightarrow[n\to\infty]{P\text{-a.s.}}0.
\end{equation*}
Combining the limits of the numerator and denominator yields
\begin{equation*}
\phi_n
\xrightarrow[n\to\infty]{P\text{-a.s.}}
\frac{
E\left[\inpr{\log_{\mu}X_1}{\log_{\mu}X_0}\right]
}{
E\left[\norm{\log_{\mu}X_0}^2\right]
}=\phi
,
\end{equation*}
which completes the proof.
\end{proof}

In the previous theorem we established the strong consistency of the estimator $\phi_n$ in the case where $\mu$ is the Fréchet mean, that is, when we have almost sure convergence to the true parameter $\mu$.
However, when $\abs{\phi}>1$, the process is no longer stationary and the previous ergodic argument does not apply.
Nevertheless, we can still obtain an estimate for $\abs{\phi}$ in certain special situations.
In particular, for manifolds with infinite injectivity radius, it is possible to estimate $\abs{\phi}$.
\begin{theorem}\label{thm: estimating absolute value of phi}
Let $M$ be a Riemannian manifold with infinite injectivity radius. 
Let $X_n$ be a R-AR$(\mu,\phi)$ with $\abs{\phi}>1$ such that $X_0$ is atomless.
Then, for all $p\in M$,
\begin{equation*}
d(X_n,p)^{1/n}\xrightarrow[n\to\infty]{P\text{-a.s.}}\abs{\phi},
\end{equation*}
as $n\to\infty$.
\end{theorem}
\begin{proof}
See Appendix \ref{sec:appendix:proofs}
\end{proof}
In the previous theorem, let $\Phi_n = d(X_n,p)^{1/n}$.
Then $\Phi_n$ provides a geometric diagnostic for nonstationarity.
Specifically, if $\Phi_n$ does not converge to a constant strictly larger than $1$, then necessarily $\abs{\phi}\leq 1$.
Thus, failure of exponential growth rules out the nonstationary regime.
Finally, although the proof exploits the linear structure of the tangent space
$T_\mu M$ via logarithmic coordinates, the conclusion is purely metric,
being expressed entirely in terms of the Riemannian distance.
This suggests that the exponential growth property may still hold for analogous models defined on more general metric spaces, such as length or geodesic spaces,
provided an appropriate notion of multiplicative dynamics is available.

\section{Numerical experiments and application}

\subsection{Simulations}
In this section, we present simulation experiments to explore the finite-sample performances of the estimators for the parameters $\mu$ and $\phi$ in a Riemannian autoregressive model and to illustrate the theoretical properties established earlier. All simulations are conducted on the two-dimensional hyperbolic space $\mathbb{H}^2$, represented via the hyperboloid model.

The data-generating process is defined as follows. 
Let $\mu \in \mathbb{H}^2$ and $\phi \in \mathbb{R}$ be fixed parameters. Let $\{\varepsilon_k\}_{k \ge 0}$ be a sequence of i.i.d.\ random vectors in the tangent space $T_\mu \mathbb{H}^2$.
Each $\varepsilon_k$ is sampled uniformly from the Euclidean disk of radius $\delta > 0$ centred at the origin of $T_\mu \mathbb{H}^2$.
To generate the noise, we construct an orthonormal basis $\{e_1, e_2\}$ of $T_\mu \mathbb{H}^2$ using the Lorentzian inner product and a Gram--Schmidt procedure. 
Each noise vector is then given by
\begin{equation*}
\varepsilon_k = R \bigl(\cos(\theta)\, e_1 + \sin(\theta)\, e_2\bigr),
\end{equation*}
where $R=\delta\sqrt{U}$, $U\sim\mathrm{Unif}(0,1)$, and $\theta\sim\mathrm{Unif}(0,2\pi)$, with $U$ and $\theta$ independent.
A trajectory $\{X_k^{(i)}\}_{k=0}^n$ is initialised as $X_0^{(i)} = \exp_\mu(\varepsilon_0^{(i)})$, and evolves according to
\begin{equation*}
X_k^{(i)}
=
\exp_\mu\left(
\phi \log_\mu(X_{k-1}^{(i)}) + \varepsilon_k^{(i)}
\right),
\qquad k = 1,\ldots,n.
\end{equation*}
For each value of $\phi$, we generate $m$ independent trajectories $\{X_k^{(i)}\}_{k=0}^n$, $i=1,\ldots,m$, using independent noise sequences.

For each trajectory, the estimators are evaluated along a sequence of increasing sample sizes $n'$. 
For each such $n'$, we compute the Fréchet mean $\hat{\mu}_{n'}^{(i)}$ based on the first $n'$ observations of the trajectory and evaluate the estimation error
\begin{equation*}
d_{\mathbb{H}^2}(\mu, \hat{\mu}_{n'}^{(i)}).
\end{equation*}

The resulting sequences of errors, plotted as functions of $n'$, illustrate the pathwise convergence behaviour of the estimators as well as their variability across independent realisations.

We now describe the estimation of the autoregressive parameter $\phi$. 
For each trajectory and each sample size $n'$, we first compute the sample Fréchet mean $\hat{\mu}_{n'}^{(i)}$ based on the first $n'$ observations. 
We then define
\begin{equation*}
\hat{\phi}_{n'}^{(i)}=\frac{\sum_{k=1}^{n'-1}
\inpr{\log_{\hat{\mu}_{n'}^{(i)}}\!\left(X_{k+1}^{(i)}\right)}{\log_{\hat{\mu}_{n'}^{(i)}}\!\left(X_k^{(i)}\right)}}{\sum_{k=1}^{n'-1}
\norm{\log_{\hat{\mu}_{n'}^{(i)}}\!\left(X_k^{(i)}\right)}^{2}}.
\end{equation*}
To assess performance, we plot the absolute value of the error, $|\hat{\phi}_{n'}^{(i)} - \phi|$, as a function of $n'$.

\begin{figure}[H]
    \centering
    \includegraphics[width=1\linewidth]{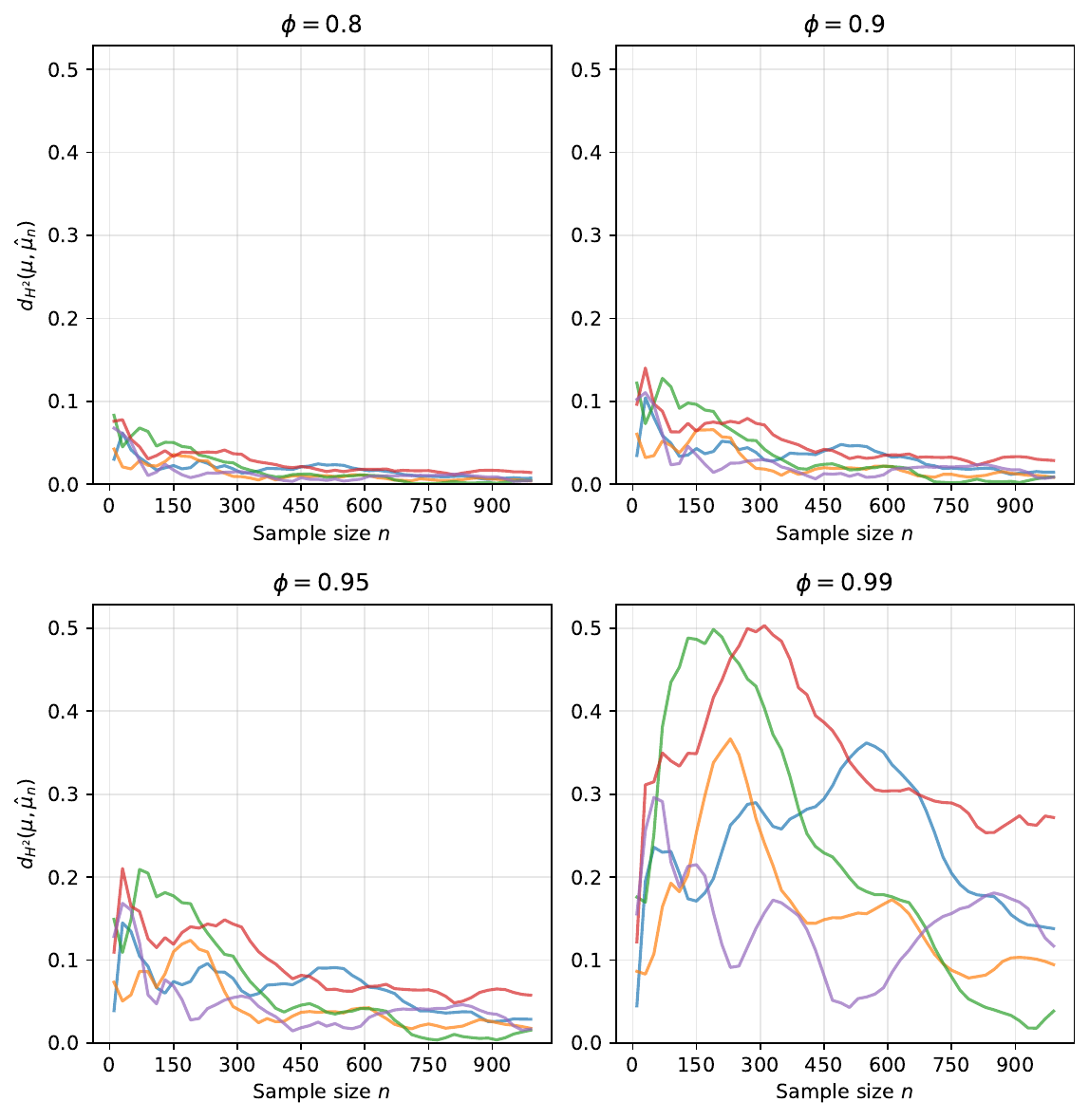}
    \caption{Estimation error $d_{\mathbb{H}^2}(\mu,\hat{\mu}_n)$ of the sample Fréchet mean as a function of the sample size $n$, for $\phi \in \{0.8,0.9,0.95,0.99\}$.
    Each curve corresponds to one independently simulated trajectory with random initialisation.}
\label{fig:mu estimation simulation}
\end{figure}

Figure~\ref{fig:mu estimation simulation} shows that the sample Fréchet mean converges toward the true parameter $\mu$, which coincides with the Fréchet mean of the stationary distribution of the process $\{X_n\}$. 
The figure displays multiple independent trajectories.
As the sample size increases, the estimator stabilises around $\mu$. The convergence is more stable for values of $\phi$ much smaller than $1$ and becomes slower and more variable as $\phi$ approaches $1$.

\begin{figure}[H]
    \centering
    \includegraphics[width=1\linewidth]{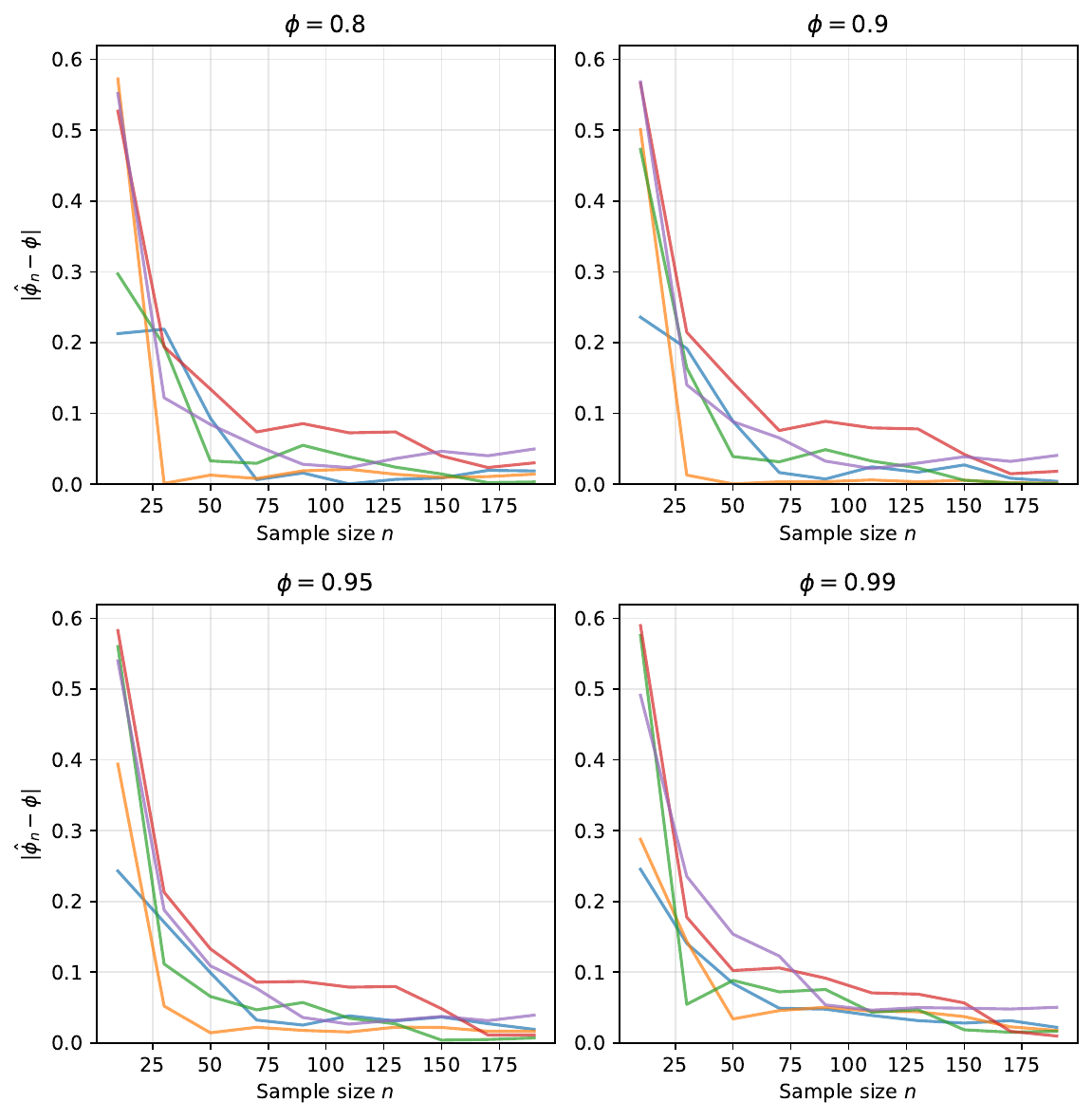}
    \caption{Absolute error $|\hat{\phi}_n-\phi|$ of the autoregressive parameter estimator as a function of the sample size $n$, for $\phi \in \{0.8,0.9,0.95,0.99\}$. 
    Each curve corresponds to one independently simulated trajectory with random initialisation.}
    \label{fig:phi estimation simulation}
\end{figure}

In Figure~\ref{fig:phi estimation simulation}, we observe the convergence of the estimator $\hat{\phi_n}$ to $\phi$ across multiple independent trajectories. 
The estimator $\hat{\phi_n}$ appears to converge to the true value faster than what we observed for the sample Fréchet mean. 
These observations are consistent with the theoretical almost sure convergence of $\hat{\phi_n}$ in the stationary case.

To summarise, the simulations support the theoretical results established for the stationary regime $|\phi|<1$. 
The sample Fréchet mean $\hat{\mu}_n$ converges toward the true parameter $\mu$, in agreement with the interpretation of $\mu$ as the Fréchet mean of the stationary law.
The convergence becomes slower and more variable as $\phi$ approaches the stationarity boundary at $1$.
The estimator $\hat{\phi}_n$ also stabilises around the true autoregressive parameter across the considered values of $\phi$.
Both finite-sample behaviours are consistent with the process remaining stationary for the considered values of $\phi$, while approaching the boundary of the stationary regime as $\phi \to 1$.

\subsection{Application}
Aerosol data analysis is fundamental in atmospheric physics, climate science, and related environmental fields; see \cite{seinfeld2016atmospheric} and \cite{hinds2022aerosol}.
In this section, we apply the R-AR model to an aerosol dataset. 
Our objective is to model the size distribution of the aerosol as an R-AR, based on sample of daily distributions of aerosol particles.
We consider the dataset on the size distribution of ambient aerosol particles collected in the UK \cite{CEDA2023}.
We focus on the 2022 year dataset, which provides continuous 5-minute measurements across the full year.
For each measurement, the independent variable is the diameter of the particles reported in discrete size bins, and the dependent variable is the corresponding aerosol number concentration, expressed as the number of particles of a given size per unit volume of air.
Measurements are provided at 5-minute intervals. 
To obtain daily data, we aggregate the particle number concentration for each diameter bin by averaging over all 5-minute measurements within a given day.
That is, we sum the intraday values and divide by the number of observations.
This yields 342 days, since some days are missing from the original dataset. 
Due to missing observations in the early part of the year, the estimation and prediction periods are taken from 1 March 2022 onward. 
Within this considered time frame, only 5 October 2022 is missing, and this date is excluded from the analysis.
Daily aerosol size distributions are often modelled using log-normal forms; see \cite{hussein2005evaluation} and \cite{finlay2020particle}. 
Although aerosol size distributions over the 14.3-679 nm range typically 
exhibit multiple modes, we adopt a 
single log-normal approximation for the daily-averaged distributions. 
This choice is motivated empirically by plotting the daily aggregated size 
distributions, we observe that they exhibit a dominant peak and are reasonably 
well described by a single log-normal curve. Accordingly, for each day we fit a 
log-normal distribution and estimate its parameters $(\mu_i, \sigma_i)$, where 
$\mu_i$ and $\sigma_i$ denote the mean and standard deviation, respectively, of 
the log-diameter distribution for day $i = 1, \ldots, 342$.

In information geometry, the parameter space of an exponential family of distributions forms a Riemannian manifold whose tensor metric, known as the Fisher--Rao metric, is given by the Fisher information matrix \cite{amari2016information}. 
For the univariate Gaussian family, this manifold is hyperbolic, and the parameter space $\set{(\mu,\sigma)}{\mu\in\mathbb{R}, \sigma>0}$ is isometric to the hyperbolic plane $\mathbb{H}^2$. 
The distance induced by the Fisher--Rao metric between two normal distributions with parameters $(\mu_1,\sigma_1)$ and $(\mu_2,\sigma_2)$ is given explicitly by
\begin{align*}
d_{\mathrm{FR}}\!\left((\mu_1,\sigma_1),(\mu_2,\sigma_2)\right)
&=\sqrt{2}\,d_{\mathbb{H}^2}\!\left(\left(\tfrac{\mu_1}{\sqrt{2}},\sigma_1\right),\left(\tfrac{\mu_2}{\sqrt{2}},\sigma_2\right)\right)
\\
&=2\sqrt{2}\operatorname{artanh}\!\left(
\sqrt{\frac{(\mu_1-\mu_2)^2 + 2(\sigma_1-\sigma_2)^2}
{(\mu_1-\mu_2)^2 + 2(\sigma_1+\sigma_2)^2}}
\right),
\end{align*}
where the hyperbolic model used here is the upper half-plane \cite{miyamoto2024closed}.
Because the Fisher--Rao metric is invariant under smooth reparametrisations, the Fisher--Rao distance between two log-normal distributions is exactly the Fisher--Rao distance between their corresponding Gaussian distributions.
Consequently, each daily aerosol size distribution can be represented as a point in $\mathbb{H}^2$, and the R-AR model can be applied to their temporal evolution as a time series on the hyperbolic plane.

We adopt three estimation windows for parameter fitting and three prediction windows for forecasting. 
The estimation window lengths are $L_{\mathrm{est}} \in \{14, 30, 60\}$ 
and the prediction window lengths are $L_{\mathrm{pred}} \in \{1, 3, 7\}$. 
We evaluate all nine combinations of the estimation and prediction windows.
Together, these choices yield the nine schemes $\mathrm{E}14\text{--}\mathrm{P}1$, $\mathrm{E}14\text{--}\mathrm{P}3$, $\ldots$, $\mathrm{E}60\text{--}\mathrm{P}7$.
To compare the forecasting performance across these nine schemes, we compute the Fisher--Rao forecast errors for each estimation--prediction pair over the full year of data.

\begin{figure}[H]
    \centering
    \includegraphics[width=\textwidth,height=0.38\textheight,keepaspectratio]{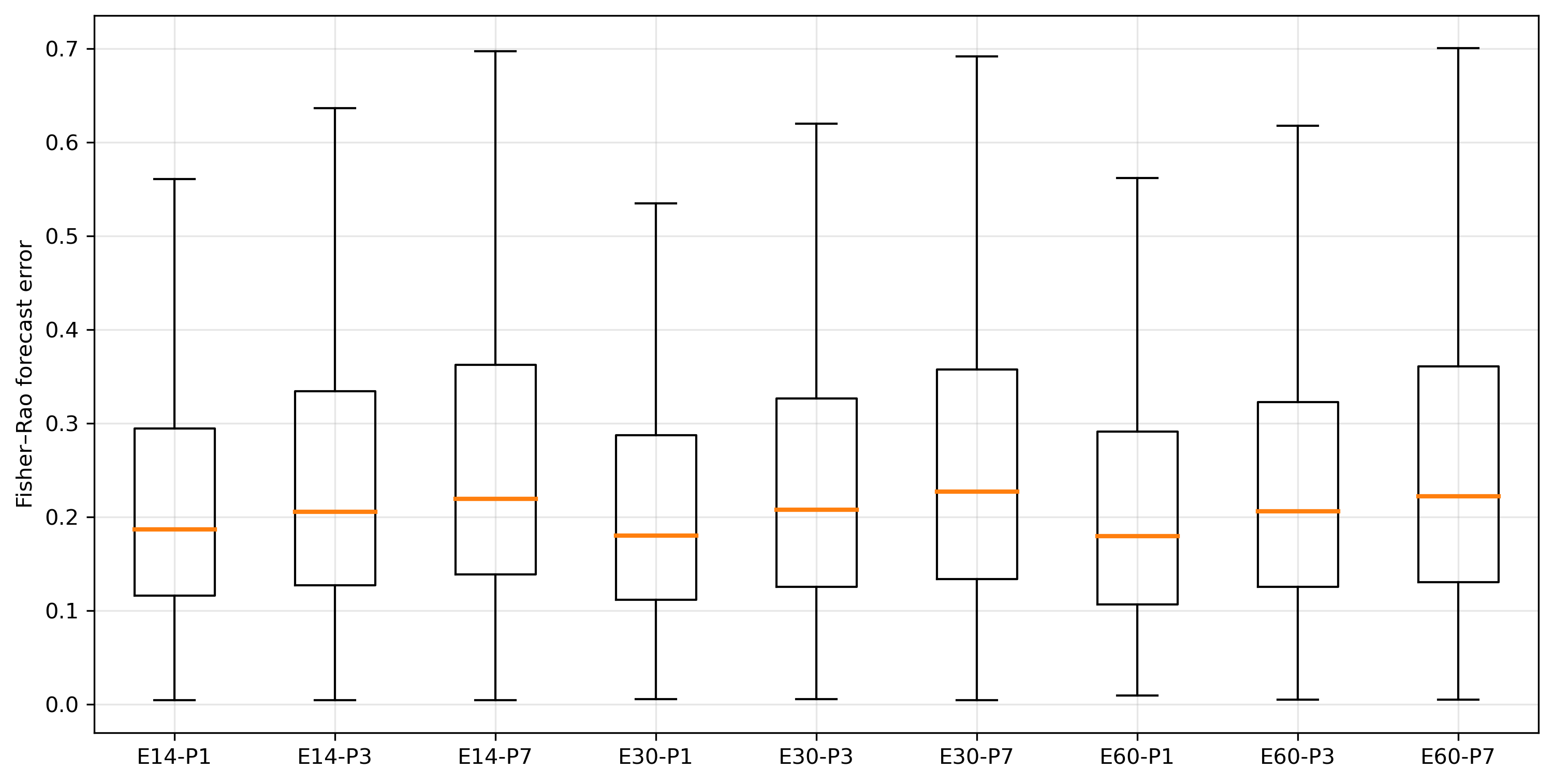}
    \caption{Boxplots of Fisher--Rao forecast errors across various estimation windows and prediction horizons. The boxes represent the interquartile range (25th to 75th percentiles), and the orange line denotes the median.}
    \label{fig:boxplot}
\end{figure}

\begin{figure}[H]
    \centering
    \includegraphics[width=\textwidth,height=0.38\textheight,keepaspectratio]{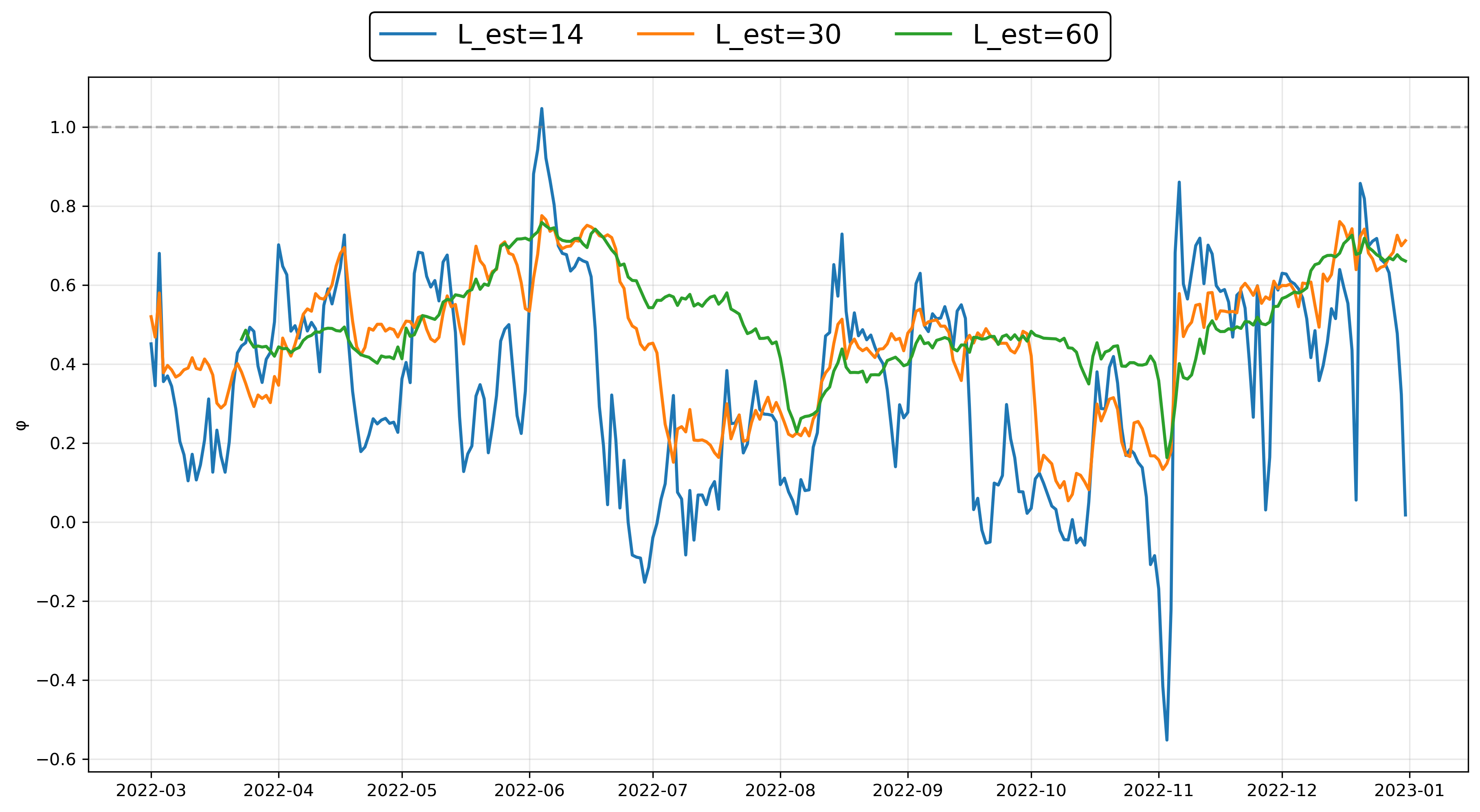}
    \caption{Estimated trajectories of the autoregressive coefficient $\phi$ across the year for the E14, E30, and E60 estimation windows.}
    \label{fig:phi_plot}
\end{figure}

\begin{figure}[H]
    \centering

    \begin{subfigure}{\textwidth}
        \includegraphics[width=\textwidth]{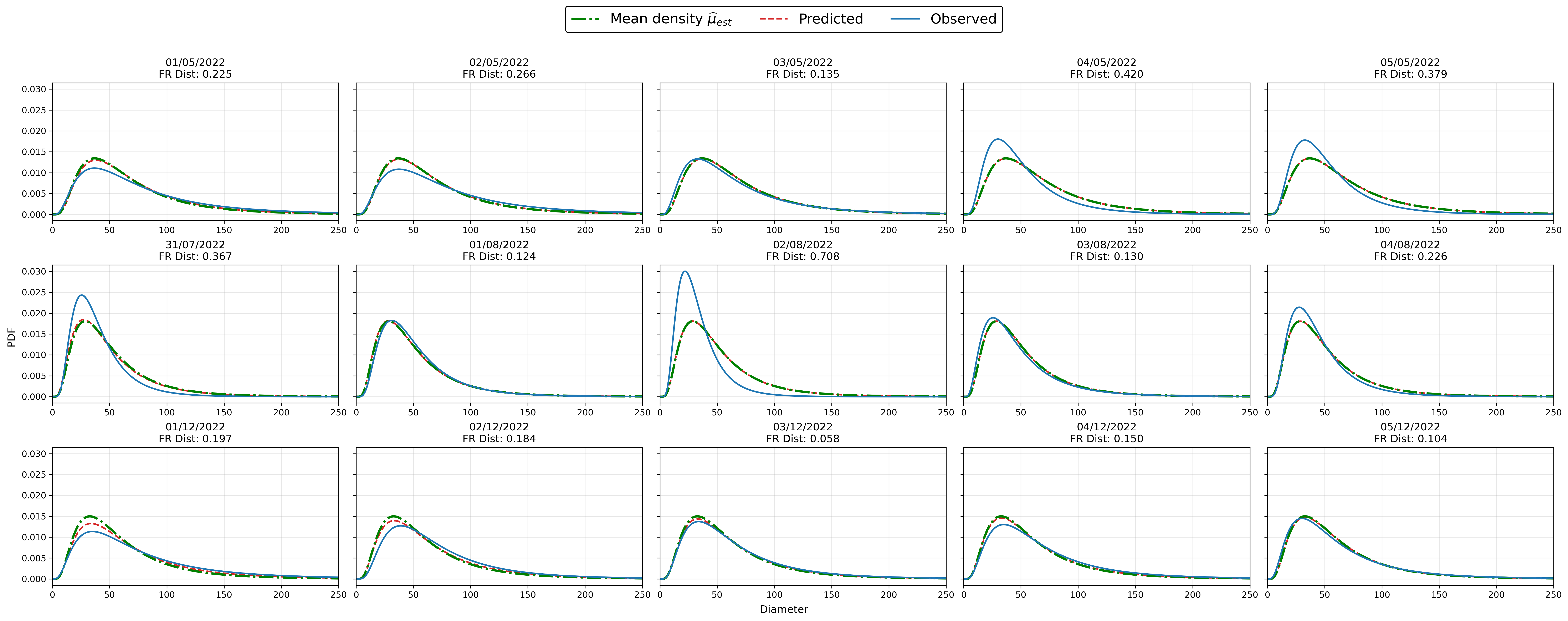}
        \caption{Comparison of estimated mean, predicted, and actual densities using a 30-day estimation window for a 5-day prediction horizon across various prediction start dates. 
        The Fisher–Rao distance for each prediction error is reported above the corresponding plot.}
        \label{fig:predicted densities}
    \end{subfigure}
    
    \vspace{1em}
    
    \begin{subfigure}{\textwidth}
        \includegraphics[width=\textwidth]{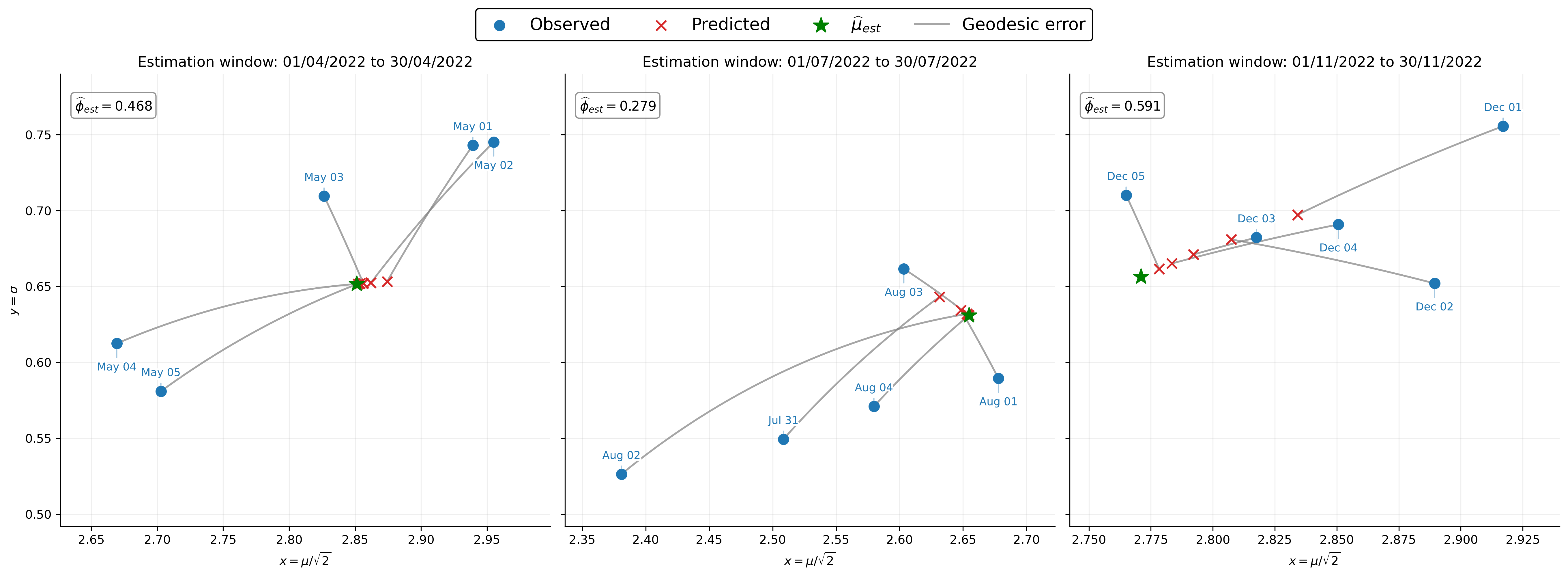}
        \caption{Predicted and observed points in the upper half-plane model of hyperbolic space.
        Each point represents a log-normal density with parameters $(\mu/\sqrt{2}, \sigma)$, and the curves indicate the hyperbolic geodesic errors between observed and predicted densities.}
        \label{fig:upper half plane prediction vs actual}
    \end{subfigure}
    
    \caption{Illustration of the E30 estimation for a five-day prediction scheme applied at various prediction start dates.
    The density plots compare predicted and observed log-normal distributions, while the upper half-plane plots show the corresponding predicted and observed points along with the geodesic prediction error.}
\label{fig: densities and upper half-plane fitting}
\end{figure}

Figure~\ref{fig:boxplot} shows the forecast error distributions using boxplots.
Forecast errors increase in both magnitude and variability as the prediction horizon lengthens, which is expected in time-series forecasting due to the accumulation of uncertainty.
Increasing the estimation window from 14 to 60 days does not lead to significant improvements in predictive accuracy.
This weak sensitivity to the window length suggests that the underlying dynamics exhibit short-term memory.
Figure~\ref{fig:phi_plot} illustrates the estimated autoregressive coefficient $\phi$ obtained using the E14, E30, and E60 estimation windows over the year.
All trajectories remain below one, indicating mean-reverting behaviour in $(\mu_t,\sigma_t)$.
The E14 estimates appear more volatile than those from E30 and E60, although all three exhibit similar temporal dynamics.
Figure~\ref{fig:predicted densities} compares selected predicted and observed densities under the E30 scheme, while Figure~\ref{fig:upper half plane prediction vs actual} displays the corresponding predicted and observed points in the upper half-plane model of hyperbolic space.
The predicted densities are seen to revert toward the mean density associated with $\widehat{\mu}_{est}$, a behaviour that is consistent with the estimated values of $\widehat{\phi}_{est}$ shown in the upper-half-plane plots.
Together, Figures~\ref{fig:predicted densities} and~\ref{fig:upper half plane prediction vs actual} show that the predicted densities are close to the estimated mean density, while the observed densities exhibit additional variability.
This indicates that the short-term prediction errors arise from variability in the daily density that is not fully captured by the R-AR model and may be better accounted for by more sophisticated models.

\vspace{-1em}

\section{Conclusions and future work}
This work presented a first-order autoregressive framework for Riemannian manifold-valued data, with the stochastic recursion formulated in a fixed tangent space and the main geometric challenges addressed in the estimation of the model parameters.
We examined the parameters involved in the model and their roles in influencing its behaviour.
We established strong consistency of the model parameter estimators using a strong law of large numbers for sample Fréchet means of ergodic Markov chains in proper metric spaces.
As numerical evidence, we provided simulations of the R-AR model in the hyperbolic plane and confirmed the validity of the estimators.
Finally, we considered aerosol data as Fisher-Rao manifold-valued dynamics, illustrating an application of the proposed model.
Overall, the proposed framework provides a geometrically coherent approach for modelling dependent data evolving on nonlinear spaces.

Several directions for future research are natural.
One direct extension is to develop a fully vector-autoregressive analogue of the present model.
In such a formulation, the autoregressive coefficient would no longer be restricted to a scalar parameter, instead it would be replaced by a linear operator acting on the tangent space at the reference point.
This would allow for richer dependence structures and could also be extended to higher-order autoregressive dynamics.
Another natural direction for future work concerns the estimation of the reference point in the expanding, nonstationary regime, where $\abs{\phi}>1$. 
In the stationary regime, the sample Fréchet mean provides a natural estimator of $\mu$, since $\mu$ can be interpreted as the Fréchet mean of the invariant distribution of the model. 
By contrast, when $\abs{\phi}>1$, the process is nonstationary, and $\mu$ no longer plays the role of the Fréchet mean of an invariant law. 
Consequently, the sample Fréchet mean is not expected to provide a consistent estimator of $\mu$ in this regime.
An alternative approach is to characterise $\mu$ through the autoregressive relation itself. 
Specifically, one may regard $\mu$ as the reference point for which the logarithmic displacements of the process satisfy the tangent-space autoregressive equation with minimal residual error. 
Under this perspective, $\mu$ is not a centre of mass, but rather the point at which the process admits the best linear autoregressive representation in $T_\mu M$. 
Developing estimators based on this characterisation, and establishing their consistency and asymptotic properties, remains an interesting problem for future investigation.

\appendix

\section{Proofs}\label{sec:appendix:proofs}

\subsection{Proof of Lemma \ref{lemma: properties of Fréchet function}}
\begin{proof}
\begin{enumerate}
\item
Assume $\mathcal{F}(p_0)<\infty$ for some $p_0\in M$.
Then, for all $p,x\in M$, using the inequality
\begin{equation*}
d^2(x,p)\leq 2d^2(x,p_0)+2d^2(p,p_0),
\end{equation*}
we obtain
\begin{equation*}
\mathcal{F}(p)=\int_{M}d^2(p,x)\lambda(dx)\leq 2\int_{M}d^2(p_0,x)\lambda(dx)+2d^2(p_0,p)<\infty.
\end{equation*}
Therefore, $\mathcal{F}(p)<\infty$ for all $p\in M$. 
\item 
Let $p,q\in M$.
\begin{align*}
\abs{\mathcal{F}(p)-\mathcal{F}(q)}&=\abs{\int_{M}\left(d^2(p,x)-d^2(q,x)\right)\lambda(dx)}\\
&\leq \int_{M}\abs{d^2(p,x)-d^2(q,x)}\lambda(dx)\\
&=\int_{M}\abs{d(p,x)-d(q,x)}(d(p,x)+d(q,x))\lambda(dx)\\
&\leq d(p,q)\int_{M}(d(x,p)+d(x,q))\lambda(dx)\\
&=d(p,q)(E[d(X,p)]+E[d(X,q)]).
\end{align*}
By Jensen’s inequality, $E[d(X,p)]^2\leq E[d^2(X,p)]<\infty$.
Moreover, the map $p\mapsto E[d(X,p)]$ is continuous on $M$ since 
$\abs{E[d(X,p)]-E[d(X,q)]}\leq d(p,q)$ for all $p\in M$.
Hence,
\begin{equation*}
\abs{\mathcal{F}(p)-\mathcal{F}(q)}\leq d(p,q)(E[d(X,p)]+E[d(X,q)])\to 0,
\end{equation*}
as $p\to q$.
\item
Let $m=\inf\set{\mathcal{F}(p)}{p\in M}$.
By continuity, $\mathcal{F}^{-1}(m)$ is closed.
Let $p,q\in \mathcal{F}^{-1}(m)$.
Then,
\begin{equation*}
d^2(p,q)\leq 2d^2(p,x)+2d^2(q,x)\implies d^2(p,q)\leq 2\mathcal{F}(p)+2\mathcal{F}(q)=4m,
\end{equation*}
so $d(p,q)\leq 2m^{1/2}$.
Thus, $\mathcal{F}^{-1}(m)$ is closed and bounded, hence compact since $M$ is proper.
To show it is nonempty, let $\{p_n\}$ be a sequence in $M$ such that $\mathcal{F}(p_n)\to m$.
Then, there exists $N\in\mathbb{N}$ such that for all $n\geq N$,
\begin{equation*}
\abs{\mathcal{F}(p_n)-m}\leq 1.
\end{equation*}
Fix $z\in M$.
Then,
\begin{equation*}
d^2(p_n,z)\leq 2d^2(p_n,x)+2d^2(x,z)
\implies d^2(p_n,z)\leq 2\mathcal{F}(p_n)+2\mathcal{F}(z)
\leq 2(m+1)+2\mathcal{F}(z).
\end{equation*}
Hence, for all $n\geq N$, $p_n\in\overline{B}_r(z)$, where 
$r^2=2(m+1)+2\mathcal{F}(z)$.
As $M$ is proper, $\overline{B}_r(z)$ is compact.
Thus, $\{p_n\}$ has a convergent subsequence $p_{n_k}\to p^*$ with $\mathcal{F}(p^*)=m$.
\end{enumerate}    
\end{proof}

\subsection{Proof of Lemma \ref{lemma: bound on the inner product of logs}}
\begin{proof}
Fix $x \in B_\mu(r^*)$. 
Define the function $f(z) = \frac{1}{2}d^2(z,x)$. 
By Theorem 6.6.1, \cite{jost2008riemannian},
\begin{equation*}
\operatorname{grad} f(z) = -\log_z x. 
\end{equation*}

Let $p,q \in S$, and let $\gamma\colon [0,\ell] \to M$ be the unique minimising unit-speed geodesic from $p$ to $q$, where $\ell = d(p,q)$. 
Define the vector field $V(s)=\log_{\gamma(s)} x$.
Since
\begin{equation*}
\log_{p}x- P_{q\to p}\log_q x=\int_{0}^\ell P_{\gamma(s)\to p}\frac{D}{ds}V(s)ds,
\end{equation*}
using the integral estimate and the fact that parallel transport is an isometry between tangent spaces,
\begin{equation*}
\norm{\log_p x - P_{q\to p}\log_q x}_p \leq \int_0^{\ell} \norm{\frac{D}{ds}V(s)}_{\gamma(s)} ds.
\end{equation*}

The integrand is the covariant derivative of the vector field $V(s) = \log_{\gamma(s)} x$ along the geodesic $\gamma$. 
By the properties of the Levi-Civita connection, this is the negative Hessian of $f$ acting on the velocity vector $\gamma'(s)$
\begin{equation*}
\frac{D}{ds}V(s) = -\nabla_{\gamma'(s)} \operatorname{grad} f = -\operatorname{Hess} f \big|_{\gamma(s)}(\gamma'(s)).
\end{equation*}

By Theorem 6.6.1 in \cite{jost2008riemannian}, since the sectional curvatures are bounded, the operator norm of the Hessian grows at most linearly with the distance to $x$. 
Specifically, there exists a constant $C'_S> 0$ depending on $S$ such that
\begin{equation*}
\norm{\frac{D}{ds}V(s)}_{\gamma(s)} \leq C'_S (1 + d(\gamma(s), x)).
\end{equation*}

Define $D_S = \max_{z \in S} d(\mu, z) < \infty$. Since $D_S \leq r^*$, the radius hypothesis guarantees the closed ball $\overline{B}_\mu(D_S)$ is strongly convex. 
Because $p, q \in S \subset \overline{B}_\mu(D_S)$, the minimising geodesic $\gamma$ is entirely contained within $\overline{B}_\mu(D_S)$. 
This ensures $d(\gamma(s), \mu) \leq D_S$ for all $s$. Using the triangle inequality $d(\gamma(s), x) \leq d(\gamma(s), \mu) + d(\mu, x) \leq D_S + d(\mu, x)$, the pointwise bound becomes
\begin{equation*}
\norm{\frac{D}{ds}V(s)}_{\gamma(s)} \leq C'_S (1 + D_S + d(\mu, x)).
\end{equation*}

Integrating this bound along the geodesic $\gamma$ of length $\ell = d(p,q)$ yields the integral estimate
\begin{equation*}
\begin{aligned}
\norm{\log_p x - P_{q\to p}\log_q x}_p &\leq \int_0^\ell C'_S (1 + D_S + d(\mu, x)) \, ds \\
&= C'_S d(p, q) (1 + D_S + d(\mu, x)) \\
&\leq C'_S (1 + D_S) d(p,q) (1 + d(\mu, x)) \\
&= C''_S d(p, q) (1 + d(\mu, x)),
\end{aligned}
\end{equation*}
where $C''_S = C'_S(1 + D_S)$ is a constant depending only on the curvature bounds and the compact set $S$.
Since parallel transport $P_{q\to p}$ is an isometry between tangent spaces, we have 
\begin{equation*}
\inpr{\log_q x}{\log_q y}_q = \inpr{P_{q\to p}\log_q x}{P_{q\to p}\log_q y}_p.
\end{equation*} 
Adding and subtracting $\inpr{P_{q\to p}\log_q x}{\log_p y}_p$, and applying the integral estimate alongside the distance bound $\norm{\log_p y}_p = d(p, y) \leq D_S + d(\mu, y)$, the inequality becomes
\begin{equation*}
\begin{aligned}
\abs{\inpr{\log_p x}{\log_p y} - \inpr{\log_q x}{\log_q y}}
&= \abs{\inpr{\log_p x}{\log_p y}
        - \inpr{P_{q\to p}\log_q x}{P_{q\to p}\log_q y}} \\
&\leq C''_S d(p,q)
\begin{aligned}[t]
\Bigl[{}& 2D_S
   + (1+D_S)\bigl(d(\mu,x)+d(\mu,y)\bigr) \\
   &+ 2d(\mu,x)d(\mu,y)
   \Bigr].
\end{aligned}
\end{aligned}
\end{equation*}

Young's inequality implies
\begin{equation*}
\begin{aligned}
\abs{\inpr{\log_p x}{\log_p y} - \inpr{\log_q x}{\log_q y}}
&\le C''_S d(p, q)
\Big[
2D_S
+ (1+D_S)\left(1 + \frac{1}{2}d^2(\mu, x)
+ \frac{1}{2}d^2(\mu, y)\right) \\
&\qquad
+ d^2(\mu, x) + d^2(\mu, y)
\Big].
\end{aligned}
\end{equation*}

Combining the constants depending only on $S$ into a constant $C_S > 0$, we obtain the final bound
\begin{equation*}
\abs{\inpr{\log_p x}{\log_p y} - \inpr{\log_q x}{\log_q y}}
\le C_S d(p, q)
\left(1 + d^2(\mu, x) + d^2(\mu, y)\right).
\end{equation*}

\end{proof}

\subsection{Proof of Theorem \ref{thm: estimating absolute value of phi}}
\begin{proof}
Let $(\Omega,\Sigma,P)$ be the underlying probability space.
Recall the model is
\begin{equation*}
X_n=\exp_{\mu}\!\bigl(\phi\,\log_{\mu}(X_{n-1})+\varepsilon_n\bigr),\qquad n\ge 1,
\end{equation*}
with i.i.d. $\varepsilon_n$ centred at $0\in T_{\mu}$ and having finite second moment. 
Assume that $\abs{\phi}>1$. 
Clearly, $\log_\mu(X_n)$ is well-defined for all $n$, since $\injr{M}=\infty$.
Iteratively, we have
\begin{equation*}
\log_{\mu}(X_n)=\phi^n\left(\log_{\mu}X_{0}+\sum_{k=1}^n \phi^{-k}\varepsilon_{k}\right).
\end{equation*}
We show that
\begin{equation*}
\log_{\mu}X_{0}+\sum_{k=1}^n \phi^{-k}\varepsilon_{k}\to S,
\end{equation*}
$P$-a.s as $n\to\infty$, for some random variable $S$.
Since $T_\mu M\simeq \mathbb R^d$ is finite-dimensional, almost sure convergence of the vector series
is equivalent to almost sure convergence of each coordinate series. 
Fix $i\in\{1,\dots,d\}$ and write
$\varepsilon_k=(\varepsilon_k^{(1)},\dots,\varepsilon_k^{(d)})$.
Define
\begin{equation*}
Y_k^{(i)}=\frac{\varepsilon_k^{(i)}}{\phi^k}.
\end{equation*}
Then $(Y_k^{(i)})_{k\in\mathbb N}$ is an independent and centred sequence of $\mathbb{R}$-valued random variables.
Since $E\|\varepsilon_1\|^2<\infty$, we have
\begin{equation*}
\sum_{k=1}^{\infty}
\var\!\left(Y_k^{(i)}\right)
=
\var(\varepsilon_1^{(i)})
\sum_{k=1}^{\infty}
|\phi|^{-2k}
<\infty,
\end{equation*}
as $|\phi|>1$.
Hence, by Kolmogorov’s convergence theorem,
\[
\sum_{k=1}^{\infty}
\frac{\varepsilon_k^{(i)}}{\phi^k}
\quad\text{converges $P$-a.s.}
\]

Since this holds for all $i=1,\dots,d$, it follows that
\begin{equation*}
\sum_{k=1}^{\infty}
\frac{\varepsilon_k}{\phi^k}
\quad\text{converges $P$-a.s. in } T_\mu M,
\end{equation*}
and thus
\begin{equation*}
\log_{\mu} X_0
+
\sum_{k=1}^n
\phi^{-k}\varepsilon_{k}
\to S
\quad\text{$P$-a.s.}
\end{equation*}
for some random vector $S \in T_\mu M$ as $n\to\infty$.
Let $Z=\sum_{k=1}^\infty \phi^{-k}\varepsilon_k$.
Since $X_0$ is atomless and independent of the noise, we have
\begin{equation*}
P(S=0)
=
P\!\left(\log_{\mu}X_0 = -Z\right)
=
E\!\left[
P\!\left(\log_{\mu}X_0 = -Z \,\Big|\, Z \right)
\right].
\end{equation*}
Conditioning on $Z=z$, $z\in T_{\mu}M$, the right-hand side becomes
\begin{equation*}
E\!\left[P\!\left(\log_{\mu}X_0 = -z\right)\right].
\end{equation*}
Since $\log_{\mu}$ is a diffeomorphism and $X_0$ is atomless,
\begin{equation*}
P\left(\log_{\mu}X_0 = -z\right)=0
\end{equation*}
for every fixed $z\in T_{\mu}M$.
Hence $P(S=0)=0$ and therefore $\norm{S}>0$ $P$-a.s.

Since $\log_{\mu}X_n=\phi^n S_n$ and $S_n\to S$ $P$-a.s., there exists $N^*$ such that $n>N^*$,
\begin{equation*}
0<\frac{1}{2}\norm{S}\leq \norm{\frac{\log_{\mu}X_{n}}{\phi^n}}\leq \frac{3}{2}\norm{S}.
\end{equation*}
This gives
\begin{equation*}
\abs{\phi}^n\frac{1}{2}\norm{S}\leq \norm{\log_{\mu}X_n}\leq \abs{\phi}^n\frac{3}{2}\norm{S}.
\end{equation*}
Since $\norm{\log_{\mu}X_n}=d(X_n,\mu)$, taking logarithms and dividing by $n$ we obtain
\begin{equation*}
\log\abs{\phi}
+
\frac{1}{n}\log\left(\frac{\norm{S}}{2}\right)
\leq
\frac{\log(d(X_n,\mu))}{n}
\leq
\log\abs{\phi}
+
\frac{1}{n}\log\left(\frac{3\norm{S}}{2}\right).
\end{equation*}
As $n\to\infty$, we obtain
\begin{equation*}
\Phi_n=d(X_n,\mu)^{1/n}\xrightarrow{\text{$P$-a.s.}} \abs{\phi},
\end{equation*}
To show this holds for any other $p$, note that
\begin{equation*}
d(X_n,\mu)-d(\mu,p)\leq d(X_n,p)<d(X_n,\mu)+d(\mu,p).
\end{equation*}
Dividing by $d(X_n,\mu)$,
\begin{equation*}
1-\frac{d(\mu,p)}{d(X_n,\mu)}\leq \frac{d(X_n,p)}{d(X_n,\mu)}\leq 1+\frac{d(\mu,p)}{d(X_n,\mu)},
\end{equation*}
and since $d(X_n,\mu)\to\infty$ and $d(X_n,\mu)^{1/n}\to\abs{\phi}$, we conclude
\begin{equation*}
d(X_n,p)^{1/n}\xrightarrow{\text{$P$-a.s.}} \abs{\phi},
\end{equation*}
as $n\to\infty$.
\end{proof}

\section{Codes and reproducibility}\label{sec:appendix:codes}
All the code used in the article can be found in the GitHub repository:
\url{https://github.com/MeshalAbuqrais/R-AR-simulations-and-aerosol-size-distributions}.

\bibliographystyle{apalike}
\bibliography{references}

@article{bulte2024autoregressive,
  title={An Autoregressive Model for Time Series of Random Objects},
  author={Bult{\'e}, Matthieu and S{\o}rensen, Helle},
  journal={arXiv preprint arXiv:2405.03778},
  year={2024}
}

@book{patrangenaru2016nonparametric,
  title={Nonparametric statistics on manifolds and their applications to object data analysis},
  author={Patrangenaru, Victor and Ellingson, Leif},
  year={2016},
  publisher={CRC Press, Taylor \& Francis Group Boca Raton}
}

@book{jost2008riemannian,
  title={Riemannian geometry and geometric analysis},
  author={Jost, J{\"u}rgen},
  volume={42005},
  year={2008},
  publisher={Springer}
}

@book{Rie-Lee,
  title={Introduction to Riemannian manifolds},
  author={Lee, John M},
  volume={2},
  year={2018},
  publisher={Springer}
}

@article{bhattacharya2003large,
  title={Large sample theory of intrinsic and extrinsic sample means on manifolds},
  author={Bhattacharya, Rabi and Patrangenaru, Vic},
  journal={The Annals of Statistics},
  volume={31},
  number={1},
  pages={1--29},
  year={2003},
  publisher={Institute of Mathematical Statistics}
}

@book{benaim2022markov,
  title={Markov Chains on Metric Spaces},
  author={Bena{\"\i}m, Michel and Hurth, Tobias},
  year={2022},
  publisher={Springer}
}

@inproceedings{frechet1948elements,
  title={Les {\'e}l{\'e}ments al{\'e}atoires de nature quelconque dans un espace distanci{\'e}},
  author={Fr{\'e}chet, Maurice},
  booktitle={Annales de l'institut Henri Poincar{\'e}},
  volume={10},
  number={4},
  pages={215--310},
  year={1948}
}

@book{hsu2002stochastic,
  title={Stochastic analysis on manifolds},
  author={Hsu, Elton P},
  number={38},
  year={2002},
  publisher={American Mathematical Soc.}
}

@book{ikeda2014stochastic,
  title={Stochastic differential equations and diffusion processes},
  author={Ikeda, Nobuyuki and Watanabe, Shinzo},
  volume={24},
  year={2014},
  publisher={Elsevier}
}

@article{fiori2014auto,
  title={AUTO-REGRESSIVE MOVING-AVERAGE DISCRETE-TIME DYNAMICAL SYSTEMS AND AUTOCORRELATION FUNCTIONS ON REAL-VALUED RIEMANNIAN MATRIX MANIFOLDS.},
  author={Fiori, Simone},
  journal={Discrete \& Continuous Dynamical Systems-Series B},
  volume={19},
  number={9},
  year={2014}
}

@article{zhu2024spherical,
  title={Spherical autoregressive models, with application to distributional and compositional time series},
  author={Zhu, Changbo and M{\"u}ller, Hans-Georg},
  journal={Journal of Econometrics},
  volume={239},
  number={2},
  pages={105389},
  year={2024},
  publisher={Elsevier}
}

@article{afsari2011riemannian,
  title={Riemannian ${L}^{p}$ center of mass: existence, uniqueness, and convexity},
  author={Afsari, Bijan},
  journal={Proceedings of the American Mathematical Society},
  volume={139},
  number={2},
  pages={655--673},
  year={2011}
}

@book{sturm2003probability,
  title={Probability measures on metric spaces of nonpositive curvature},
  author={Sturm, Karl-Theodor},
  year={2003},
  publisher={SFB 611}
}

@article{pennec2006intrinsic,
  title={Intrinsic statistics on Riemannian manifolds: Basic tools for geometric measurements},
  author={Pennec, Xavier},
  journal={Journal of Mathematical Imaging and Vision},
  volume={25},
  number={1},
  pages={127--154},
  year={2006},
  publisher={Springer}
}

@article{finlay2020particle,
  title={Particle size distributions},
  author={Finlay, Warren H and Darquenne, Chantal},
  journal={Journal of aerosol medicine and pulmonary drug delivery},
  volume={33},
  number={4},
  pages={178--180},
  year={2020},
  publisher={Mary Ann Liebert, Inc., publishers 140 Huguenot Street, 3rd Floor New~…}
}

@article{hussein2005evaluation,
  title={Evaluation of an automatic algorithm for fitting the particle number size distributions},
  author={Hussein, Tareq and Dal Maso, Miikka and Pet{\"a}j{\"a}, Tuukka and Koponen, Ismo K and Paatero, Pentti and Aalto, Pasi P and H{\"a}meri, Kaarle and Kulmala, Markku},
  journal={Boreal environment research},
  volume={10},
  number={5},
  pages={337},
  year={2005},
  publisher={Finnish Environment Institute}
}

@article{miyamoto2024closed,
  title={On closed-form expressions for the Fisher--Rao distance},
  author={Miyamoto, Henrique K and Meneghetti, F{\'a}bio CC and Pinele, Julianna and Costa, Sueli IR},
  journal={Information Geometry},
  volume={7},
  number={2},
  pages={311--354},
  year={2024},
  publisher={Springer}
}

@inproceedings{lou2020differentiating,
  title={Differentiating through the fr{\'e}chet mean},
  author={Lou, Aaron and Katsman, Isay and Jiang, Qingxuan and Belongie, Serge and Lim, Ser-Nam and De Sa, Christopher},
  booktitle={International conference on machine learning},
  pages={6393--6403},
  year={2020},
  organization={PMLR}
}

@book{amari2016information,
  title={Information geometry and its applications},
  author={Amari, Shun-ichi},
  volume={194},
  year={2016},
  publisher={Springer}
}

@book{hinds2022aerosol,
  title={Aerosol technology: properties, behavior, and measurement of airborne particles},
  author={Hinds, William C and Zhu, Yifang},
  year={2022},
  publisher={John Wiley \& Sons}
}

@book{seinfeld2016atmospheric,
  title={Atmospheric chemistry and physics: from air pollution to climate change},
  author={Seinfeld, John H and Pandis, Spyros N},
  year={2016},
  publisher={John Wiley \& Sons}
}

@article{barbati1994hausdorff,
  title={The Hausdorff metric topology, the Attouch-Wets topology and the measurability of set-valued functions},
  author={Barbati, Alberto and Beer, Gerald and Hess, Christian},
  journal={J. Convex Anal},
  volume={1},
  number={1},
  pages={107--119},
  year={1994}
}

@article{evans2024limit,
  title={Limit theorems for Fr{\'e}chet mean sets},
  author={Evans, Steven N and Jaffe, Adam Q},
  journal={Bernoulli},
  volume={30},
  number={1},
  pages={419--447},
  year={2024},
  publisher={Bernoulli Society for Mathematical Statistics and Probability}
}

@article{ding2025manifold,
  title={Manifold-valued models for analysis of EEG time series data},
  author={Ding, Tao and Nye, Tom MW and Wang, Yujiang},
  journal={Computational Statistics \& Data Analysis},
  volume={209},
  pages={108168},
  year={2025},
  publisher={Elsevier}
}

@inproceedings{vemulapalli2016rolling,
  title={Rolling rotations for recognizing human actions from 3d skeletal data},
  author={Vemulapalli, Raviteja and Chellapa, Rama},
  booktitle={Proceedings of the IEEE conference on computer vision and pattern recognition},
  pages={4471--4479},
  year={2016}
}

@inproceedings{nava2025ridge,
  title={Ridge Regression for Manifold-Valued Time-Series with Application to Hurricane Forecasting},
  author={Nava-Yazdani, Esfandiar},
  booktitle={International Conference on Geometric Science of Information},
  pages={3--11},
  year={2025},
  organization={Springer}
}

@book{jelinek2017ecg,
  title={ECG time series variability analysis: engineering and medicine},
  author={Jelinek, Herbert F and Cornforth, David J and Khandoker, Ahsan H},
  year={2017},
  publisher={CRC Press}
}

@book{cipra2020time,
  title={Time series in economics and finance},
  author={Cipra, Tomas and others},
  year={2020},
  publisher={Springer}
}

@misc{CEDA2023,
  author = {{CEDA}},
  title  = {{\href{https://catalogue.ceda.ac.uk/uuid/6f6366eb78c44875b2c92d3b8fe403c1}{SMPS aerosol particle data at the Manchester Air Quality Site}}},
  year   = {2019-2023}
}

@article{lee2025kolmogorov,
  title={On the Kolmogorov-Feller weak law of large numbers for Frechet mean on non-compact symmetric spaces},
  author={Lee, Jongmin and Jung, Sungkyu},
  journal={arXiv preprint arXiv:2509.02074},
  year={2025}
}

@article{jaffe2024fr,
  title={Fr{\'e}chet Means in Infinite Dimensions},
  author={Jaffe, Adam Quinn},
  journal={arXiv preprint arXiv:2410.17214},
  year={2024}
}

\end{document}